\documentclass{bcp}
\usepackage{amsmath,amsthm}
\usepackage{amssymb}

\usepackage{graphicx}

\newtheorem{theorem}{Theorem}[section]

\newtheorem{lemma}[theorem]{Lemma}

\theoremstyle{definition}
\newtheorem{definition}[theorem]{Definition}

\newtheorem{remark}[theorem]{Remark}

\def\E{{\mathrm{E}}}

\def \R{{{\rm I{\!}\rm R}}}
\def \N{{\mathbb{N}}}
\def \Z{{\mathbb{Z}}}
\newcommand{\cc}{{\mathbb C}}

\usepackage[numbers]{natbib}
\usepackage{multicol}
\usepackage{enumerate}
\usepackage{amsfonts}
\usepackage{epstopdf}
\usepackage{setspace}
\usepackage{color}
\usepackage{xcolor}
\usepackage[toc,page]{appendix}
\usepackage{array}
\newcolumntype{C}[1]{>{\centering\arraybackslash}m{#1}}

\usepackage{diagbox}
\usepackage{makecell}
\usepackage{slashbox,pict2e}

\allowdisplaybreaks 

\begin{document}

\mathclass{Primary $62M10$; Secondary $37M10$.}

\abbrevauthors{A. Grzesiek and A. Wy{\l}oma{\'n}ska}
\abbrevtitle{Asymptotics of dependence measures for bivariate $\alpha-$stable AR(1) model}

\title{Asymptotic behavior of the cross-dependence measures for bidimensional AR(1) model with $\alpha-$stable noise}

\author{Aleksandra Grzesiek}
\address{Faculty of Pure and Applied Mathematics, Hugo Steinhaus Center\\
Wroc{\l}aw University of Science and  Technology,\\
Wybrze\.ze Wyspia{\'n}skiego 27, 50-370 Wroc{\l}aw, Poland\\
E-mail: aleksandra.grzesiek@pwr.edu.pl}

\author{Agnieszka Wy{\l}oma{\'n}ska}
\address{Faculty of Pure and Applied Mathematics, Hugo Steinhaus Center\\
	Wroc{\l}aw University of Science and  Technology,\\
	Wybrze\.ze Wyspia{\'n}skiego 27, 50-370 Wroc{\l}aw, Poland\\
	E-mail: agnieszka.wylomanska@pwr.edu.pl}

\maketitlebcp

\begin{abstract}
	In this paper, we consider a bidimensional autoregressive model of order $1$ with $\alpha-$stable noise. Since in this case the classical measure of dependence known as the covariance function is not defined, the spatio-temporal dependence structure is described using the alternative measures, namely the codifference and the covariation functions. Here, we investigate the asymptotic relation between these two dependence measures applied to the description of the cross-dependence of the bidimensional model. We demonstrate the case when the dependence measures are asymptotically proportional with the coefficient of proportionality equal to the parameter $\alpha$. The theoretical results are supported by illustrating the asymptotic behavior of the dependence measures for two exemplary bidimensional $\alpha-$stable AR(1) systems.
\end{abstract}

\section{Introduction} 
A lot of research has shown that the models based on the heavy-tailed distributions are suitable for describing {various} kinds of phenomena where the assumption about the Gaussian distribution is not reasonable, see \cite{mit,fin2,phy2,Gardes2010,Kysely2010,valueatrisk,earthquake,infinite11,panor2}. Due to the generalized central limit theorem, the natural extension of the Gaussian distribution is the $\alpha-$stable one \cite{Taqqu,zolotarev,stab3}. The $\alpha-$stable distribution (called also stable distribution) belongs to the heavy-tailed class of distributions (except the Gaussian case of $\alpha=2$) and it is characterized by infinite variance. In the literature, one can find many applications of stochastic processes and time series based on the $\alpha-$stable distribution in both one-dimensional and multidimensional case, see \cite{stab4,panor1,stab_stoyanov,stab_obuchowski,multiarmastable,ar,ARMA_measures,dia2}. In the multidimensional case, the classical examples are the financial data (which mostly are non-Gaussian) related to different assets (like exchange rates and metal prices) which are in some sense connected however their relationship is observed {with the time shift} \cite{stab_jablonska,aw2,aw3}.

However, the classical dependence measure known as the covariance function is not defined for $\alpha-$stable random vectors and thus it cannot be used to describe the structure of dependence for the $\alpha-$stable distribution-based processes.  Therefore, alternative dependence measures are applied to such models, like for example covariation (and auto-covariation) \cite{Taqqu,ar,est2}, codifference (and auto-codifference) \cite{Taqqu,ar,cod1,cod2,Chechkin2015} or fractional lower order covariance \cite{floc1,floc2,floc3}, see also \cite{spectral_cov}. The mentioned measures can be also used to describe the spatio-temporal dependence structure for multidimensional models. In the author's previous works the cross-codifference and the cross-covariation functions were introduced for the bidimensional autoregressive model \mbox{of order 1} with $\alpha-$stable noise (bidimensional $\alpha-$stable AR(1) model) to describe the dependence between the components of two-dimensional model, see \cite{nasza2,nasza}.

In the literature, one can find the research papers where the asymptotic behavior of the dependence measures is studied. In particular, the asymptotic relation between the auto-codifference and the auto-covariation functions of a given $\alpha-$stable distribution-based process is examined, see for example \cite{ar,parma,arma}. In this paper, we continue this research by studying the asymptotic relationship between the {cross-codifference} and the {cross-covariation} functions for the bidimensional $\alpha-$stable autoregressive model \mbox{of order 1}. The work can be seen as an extension of the results presented in \cite{arma}, where the asymptotic behavior of the ratio of the auto-codifference function and the auto-covariation function for one-dimensional autoregressive time series was investigated. Here we {extend} the approach to the two-dimensional model {by studying the asymptotic relation of the cross-dependence measures} describing the relationship between the components of the bidimensional model. We demonstrate the {case} when the dependence measures are asymptotically proportional with the coefficient of proportionality equal to the parameter $\alpha$. This result can be a starting point for the introduction of a new estimation method of the stability index in the considered model. 

The paper is organized as follows. In Section \ref{sec2} we present the definition of the bidimensional autoregressive model of order $1$ with the $\alpha-$stable noise. In Section \ref{measures} we recall the definitions and properties of the codifference and the covariation functions together with the formulas for the {cross-dependence} measures corresponding to the bidimensional $\alpha-$stable AR(1) model. Section \ref{section_proof} contains the theorem concerning the asymptotic behavior of the ratio of the dependence measures for the bidimensional AR(1) model with the $\alpha-$stable noise. In Section \ref{example_sec} we illustrate the theoretical results by considering two exemplary bidimensional $\alpha-$stable AR(1) systems. In Section \ref{Conclusions} we conclude the paper.

	\section{Bidimensional AR(1) model with $\alpha-$stable noise} \label{sec2}
%In this paper we study the bidimensional autoregressive model of order $1$ with $\alpha-$stable noise.
In this section, we introduce the bidimensional autoregressive model of order $1$ with $\alpha-$stable noise. Let us begin by recalling the definition of the symmetric $\alpha-$stable random vectors which are the multidimensional version of the symmetric $\alpha-$stable random variables. We remind that the $\alpha-$stable distribution with $0<\alpha\leq 2$ can be considered as the extension of the Gaussian distribution (the case of $\alpha=2$). Let $S_{d}=\{\mathbf{s}:||s||=1\}$ be a unit sphere in $\R^d$. Then a symmetric $\alpha-$stable random vector can be presented via the characteristic function.
\begin{theorem}
	\textup{\citep{Taqqu}} The vector $\mathbf{Z}=(Z_1,Z_2,\ldots,Z_d)$ is a symmetric $\alpha-$stable vector in $\R^d$ with $0<\alpha<2$ if and only if there exists a unique symmetric finite spectral measure $\Gamma(\cdot)$ on the unit sphere $S_d$ such that 
	{\small \begin{eqnarray}\label{charact2}
	\phi_{\mathbf{Z}}(\mathbf{\boldsymbol{\theta}})=\E[\exp\{i\langle \boldsymbol{\theta}, \mathbf{Z} \rangle\}]=
	\exp\left\{-\int_{S_d}{|\langle \boldsymbol{\theta},\mathbf{s} \rangle|}^{\alpha}\Gamma(ds)\right\},
	\end{eqnarray}}
	where $\langle \cdot,\cdot \rangle$ is the inner product.
	\label{theo}
\end{theorem} 
As one can see the information about the distribution is included in the spectral measure $\Gamma(\cdot)$ and the parameter $\alpha$, also called the stability index. They fully describe $d$-dimensional symmetric $\alpha-$stable distribution which is denoted as
{\small\begin{equation}
\mathbf{Z} \sim S_{\alpha}(\Gamma).
\end{equation}}
For more information about the multidimensional $\alpha-$stable distribution see for instance \cite{Taqqu}. The following definition applies to the bidimensional AR(1) time series based on the $\alpha-$stable random vector.
\begin{definition}
	The time series $\{\mathbf{X}(t)\}=\{(X_1(t),X_2(t))\}$ is \textit{a bidimensional autoregressive model of \mbox{order $1$} with $\alpha-$stable noise} if for every $t\in \Z$ it satisfies the following equation
	{\small\begin{align}
	\mathbf{X}\left(t\right)-\Theta\ \mathbf{X}\left(t-1\right)=\mathbf{Z}\left(t\right),
	\label{eq1}
	\end{align}}
	where $\Theta$ is a $2 \times 2$ coefficients' matrix given by
	{\small\begin{equation}
	\Theta= \left[
	\begin{array}{cc}
	a_1 & a_2\\
	a_3 & a_4
	\end{array}
	\right],
	\label{theta}
	\end{equation}}
	%	and $\{\mathbf{Z}(t)\}=\{Z_1(t),Z_2(t)\}$ is a bidimensional stable vector in $\R^2$ with the characteristic function defined in (\ref{char_fun}) under the assumption that $\mathbf{Z}(t)$ is independent from $\mathbf{Z}(t+h)$ for all $h \neq 0$.
	$\mathbf{Z}(t)$ is independent from $\mathbf{Z}(t+h)$ for all $h \neq 0$ and $\{\mathbf{Z}(t)\}=\{(Z_1(t),Z_2(t))\}$ is a bidimensional symmetric $\alpha-$stable vector in $\R^2$ with the characteristic function defined in Eq. (\ref{charact2}).
\end{definition}
Moreover, we assume that for the system given by Eq. (\ref{eq1}) the following condition is satisfied
{\small\begin{equation}
\text{det}(I-z\Theta)\neq0 \quad \text{for all}\ z \in \cc\ \text{such that}\ |z| \leq 1,
\label{condition}
\end{equation}}
i.e. the eigenvalues of the matrix $\Theta$ are less than $1$ in the absolute value. Under this assumption, for each $t$ we can express $\mathbf{X}(t)$ in the causal representation as
{\small\begin{align}
\mathbf{X}(t)=\sum_{j=0}^{+\infty}\Theta^j\mathbf{Z}\left(t-j\right),
\label{eq2}
\end{align}}
where the coefficients $\Theta^j$ are absolutely summable. In this case, $\mathbf{X}(t)$ is bounded (in the sense of the so-called covariation norm in the space of the $\alpha-$stable random variables, \cite{Taqqu}). In this paper, we consider only the case when the eigenvalues of the matrix $\Theta$ are real numbers.

\noindent Moreover, in \cite{matrix} it was shown that for the $2 \times 2$ matrix the coefficients of $\Theta^j$ can be expressed as 
{\small\begin{equation}
\Theta^j= \left[
\begin{array}{cc}
\frac{\lambda_2 \lambda_1^j-\lambda_1 \lambda_2^j}{\lambda_2-\lambda_1}+\frac{\lambda_2^j-\lambda_1^j}{\lambda_2-\lambda_1}a_1 & \frac{\lambda_2^j-\lambda_1^j}{\lambda_2-\lambda_1}a_2\\
\frac{\lambda_2^j-\lambda_1^j}{\lambda_2-\lambda_1}a_3 & \frac{\lambda_2 \lambda_1^j-\lambda_1 \lambda_2^j}{\lambda_2-\lambda_1}+\frac{\lambda_2^j-\lambda_1^j}{\lambda_2-\lambda_1}a_4
\end{array}
\right],
\label{theta1}
\end{equation}}
where $\lambda_1$, $\lambda_2$ are two different eigenvalues of the matrix $\Theta$, i.e. when $(a_1-a_4)^2 > -4a_2a_3$ (and $|\lambda_1|<1$, $|\lambda_2|<1$) or 
{\small\begin{equation}
\Theta^j= \left[
\begin{array}{cc}
j\lambda^{j-1}a_1-(j-1)\lambda^j & j\lambda^{j-1}a_2\\
j\lambda^{j-1}a_3 & j\lambda^{j-1}a_4-(j-1)\lambda^j
\end{array}
\right],
\label{theta2}
\end{equation}}
where the eigenvalues of the matrix $\Theta$ are equal $\lambda_1=\lambda_2=\lambda$, i.e. when $(a_1-a_4)^2 = -4a_2a_3$ (and $|\lambda|<1$).
	\section{Measures of dependence for bidimensional AR(1) model with $\alpha-$stable noise} \label{measures}
In the case of the $\alpha-$stable random vectors (for $\alpha<2$), the classical dependence measure known as the covariance function is not defined due to the infinite second moment and therefore other measures of dependence have to be used. The most popular ones are the codifference and the covariation given in Definition \ref{def_cod} and Definition \ref{def_cov}, respectively.

\begin{definition} \label{def_cod}
	%Let the random vector $(Z_1,Z_2)$ be a bidimensional $\alpha-$stable vector with \mbox{$0 <\alpha < 2$} and the spectral measure $\Gamma(\cdot)$. Then the codifference between $Z_1$ and $Z_2$ is given by \textup{\citep{Taqqu,Chechkin2015}}
	\textup{\citep{Taqqu,Chechkin2015}} Let us consider the random vector $(Z_1,Z_2)$. Then the \textit{codifference} between $Z_1$ and $Z_2$ is given by
	{\small\begin{equation}
	\mathrm{CD}(Z_1,Z_2)=\log \E\exp\{i(Z_1-Z_2)\}-\log \E\exp\{iZ_1\} -\log \E\exp\{-iZ_2\}.\label{coddiference2}
	\end{equation}}
\end{definition}

\begin{definition} \label{def_cov}
	\textup{\citep{Taqqu}} Let us consider the bidimensional symmetric $\alpha-$stable random vector $(Z_1,Z_2)$ with \mbox{$1 <\alpha < 2$} and let $\Gamma(\cdot)$ be the spectral measure of $(Z_1,Z_2)$. The \textit{covariation} of $Z_1$ on $Z_2$ is the real number defined as 
	{\small\begin{equation}
	\mathrm{CV}(Z_1,Z_2)=\int_{S_2}s_1s_2^{\langle \alpha-1 \rangle}\Gamma(ds), \label{covariation}
	\end{equation}}
	where $a^{\langle p\rangle}$ is called the signed power and it is equal to {\small$$a^{\langle p\rangle}=|a|^p\rm{sign}(a).$$}
\end{definition}

Let us note that the covariation is defined only for the symmetric $\alpha-$stable random vectors. Moreover, the covariation function is not symmetric in its arguments, in contrast to the codifference function which is symmetric for the symmetric random vectors. It is worth mentioning, that in the case of two independent random variables $Z_1$ and $Z_2$ both measures are equal to $0$, i.e. $\mathrm{CD}(Z_1,Z_2)=\mathrm{CV}(Z_1,Z_2)=0$. Moreover, for the Gaussian random vectors $(Z_1,Z_2)$ both measures reduce to the classical covariance function, namely {\small$$\mathrm{Cov}(Z_1,Z_2)=\mathrm{CD}(Z_1,Z_2)=2\mathrm{CV}(Z_1,Z_2),$$}
where $\mathrm{Cov}(Z_1,Z_2)$ denotes the covariance function. More properties of the codifference and the covariation one can find in \cite{Taqqu}.

The functions defined in Eq. (\ref{coddiference2}) and in Eq. (\ref{covariation}) can be used to describe the interdependence of a stochastic process $\{X(t)\}$ as the \mbox{auto-codifference} or the auto-covariation, see \cite{arma,parma,ar,cod2,dep7,Chechkin2015}. In the authors' previous papers the measures were also applied to describe the spatio-temporal dependence structure of the bidimensional AR(1) model $\{(X_1(t),X_2(t))\}$ as the cross-codifference and the cross-covariation, see \cite{nasza2,nasza}.

Now, for the bidimensional AR(1) model with $\alpha-$stable noise defined in Eq. (\ref{eq2}) we can rewrite the formulas for the \mbox{cross-dependence} measures given in \cite{nasza2} taking into account the expressions for the $j$-th power of the coefficients' matrix $\Theta$ given in Eq. (\ref{theta1}) or Eq. (\ref{theta2}). {The formulas for the cross-codifference and the cross-covariation are presented in Lemma \ref{th1} and Lemma \ref{th2}, respectively.}

{\begin{lemma}
		\label{th1} 	Let $\{\mathbf{X}(t)\}=\{X_1(t),X_2(t)\}$ be the bounded solution of Eq. (\ref{eq1}) given by Eq. (\ref{eq2}). Let $t\in \Z$ and $h \in \N_{0}=\N\cup \{0\}$, then
		\begin{enumerate}
			\item in the case of two different eigenvalues of the coefficients' matrix $\Theta$ denoted as $\lambda_1$ and $\lambda_2$, where $|\lambda_1|<1$, $|\lambda_2|<1$, and $0<\alpha<2$ we have 
			{\small%\begin{equation} 
			\begin{align} 
			\rm{CD}&(X_1(t),X_2(t-h))\nonumber\\&=\sum_{j=0}^{+\infty}\int_{S_2}\Bigg(\left|{\frac{\lambda_1^h(\lambda_2\lambda_1^js_1-\lambda_1^ja_1s_1-\lambda_1^ja_2s_2)+\lambda_2^h(-\lambda_1\lambda_2^js_1+\lambda_2^ja_1s_1+\lambda_2^ja_2s_2)}{\lambda_2-\lambda_1}}\right|^{\alpha}\nonumber\\&+\left|\frac{\lambda_2^ja_3s_1-\lambda_1^ja_3s_1+\lambda_2\lambda_1^js_2-\lambda_1\lambda_2^js_2+\lambda_2^ja_4s_2-\lambda_1^ja_4s_2}{\lambda_2-\lambda_1}\right|^{\alpha}\label{eq3} \\&-\Bigg|\frac{\lambda_2^ja_3s_1-\lambda_1^ja_3s_1+\lambda_2\lambda_1^js_2-\lambda_1\lambda_2^js_2+\lambda_2^ja_4s_2-\lambda_1^ja_4s_2}{\lambda_2-\lambda_1}\nonumber\\&+\frac{\lambda_1^h(-\lambda_2\lambda_1^js_1+\lambda_1^ja_1s_1+\lambda_1^ja_2s_2)+\lambda_2^h(\lambda_1\lambda_2^js_1-\lambda_2^ja_1s_1-\lambda_2^ja_2s_2)}{\lambda_2-\lambda_1}\Bigg|^\alpha\Bigg)\Gamma(ds) \nonumber
			\end{align}}
			%\end{equation}}
			{\small%\begin{equation}
			\begin{align} 
			\rm{CD}&(X_1(t),X_2(t+h))\nonumber\\&=\sum_{j=0}^{+\infty}\int_{S_2}\Bigg(\left|\frac{\lambda_1^h(-\lambda_1^ja_3s_1+\lambda_2\lambda_1^js_2-\lambda_1^ja_4s_2)+\lambda_2^h(\lambda_2^ja_3s_1-\lambda_1\lambda_2^js_2+\lambda_2^ja_4s_2)}{\lambda_2-\lambda_1}\right|^{\alpha}\nonumber\\&+\left|\frac{\lambda_2\lambda_1^js_1-\lambda_1\lambda_2^js_1+\lambda_2^ja_1s_1-\lambda_1^ja_1s_1+\lambda_2^ja_2s_2-\lambda_1^ja_2s_2}{\lambda_2-\lambda_1}\right|^{\alpha}\label{eq7}\\&-\Bigg|\frac{\lambda_2\lambda_1^js_1-\lambda_1\lambda_2^js_1+\lambda_2^ja_1s_1-\lambda_1^ja_1s_1+\lambda_2^ja_2s_2-\lambda_1^ja_2s_2}{\lambda_2-\lambda_1}\nonumber
			\\&+\frac{\lambda_1^h(\lambda_1^ja_3s_1-\lambda_2\lambda_1^js_2+\lambda_1^ja_4s_2)+\lambda_2^h(-\lambda_2^ja_3s_1+\lambda_1\lambda_2^js_2-\lambda_2^ja_4s_2)}{\lambda_2-\lambda_1}\Bigg|^\alpha\Bigg)\Gamma(ds)\nonumber
			\end{align}}
		%	\end{equation}
			\item in the case of equal eigenvalues of the coefficients' matrix $\Theta$ denoted as $\lambda_1=\lambda_2=\lambda$, where $|\lambda|<1$, and $0<\alpha<2$ we have 
			{\small%\begin{equation}
			\begin{align} 
			{\rm{CD}}&(X_1(t),X_2(t-h))=\sum_{j=0}^{+\infty}\int_{S_2}(|\lambda^h{(j\lambda^{j-1}a_1s_1-(j-1)\lambda^js_1+j\lambda^{j-1}a_2s_2)}\label{eq5}\\&+h\lambda^h{(\lambda^{j-1}a_1s_1-\lambda^js_1+\lambda^{j-1}a_2s_2)}|^{\alpha}+|j\lambda^{j-1}a_3s_1-(j-1)\lambda^js_2+j\lambda^{j-1}a_4s_2|^{\alpha}\nonumber\\&-|j\lambda^{j-1}a_3s_1-(j-1)\lambda^js_2+j\lambda^{j-1}a_4s_2-\lambda^h{(j\lambda^{j-1}a_1s_1-(j-1)\lambda^js_1+j\lambda^{j-1}a_2s_2)}\nonumber\\&-h\lambda^h{(\lambda^{j-1}a_1s_1-\lambda^js_1+\lambda^{j-1}a_2s_2)}|^\alpha)\Gamma(ds). \nonumber
			\end{align}}
		%	\end{equation}
			{\small%\begin{equation}
			\begin{align} 
			{\rm{CD}}&(X_1(t),X_2(t+h)=\sum_{j=0}^{+\infty}\int_{S_2}(|\lambda^h(j\lambda^{j-1}a_3s_1-(j-1)\lambda^js_2+j\lambda^{j-1}a_4s_2))\label{eq9}\\&+h\lambda^h(\lambda^{j-1}a_3s_1+\lambda^{j-1}a_4s_2-\lambda^{j}s_2)|^{\alpha}+|j\lambda^{j-1}a_1s_1-(j-1)\lambda^{j}s_1+j\lambda^{j-1}a_2s_2|^{\alpha}\nonumber\\&-|j\lambda^{j-1}a_1s_1-(j-1)\lambda^{j}s_1+j\lambda^{j-1}a_2s_2-\lambda^h(j\lambda^{j-1}a_3s_1-(j-1)\lambda^js_2+j\lambda^{j-1}a_4s_2)\nonumber\\&-h\lambda^h(\lambda^{j-1}a_3s_1+\lambda^{j-1}a_4s_2-\lambda^{j}s_2)|^\alpha)\Gamma(ds)\nonumber
			\end{align}}
			%\end{equation}	
		\end{enumerate}
\end{lemma}}

\begin{proof}
	The above equations follow from the formulas presented in the authors' previous paper, see \cite{nasza2}, and the expressions for the $j$-th power of the coefficients' matrix $\Theta$ given in Eq. (\ref{theta1}) or Eq. (\ref{theta2}).
\end{proof}

{\begin{lemma}
		\label{th2} Let $\{\mathbf{X}(t)\}=\{X_1(t),X_2(t)\}$ be the bounded solution of Eq. (\ref{eq1}) given by Eq. (\ref{eq2}). Let $t\in \Z$ and $h \in \N_{0}=\N\cup \{0\}$, then
		\begin{enumerate}
			\item in the case of two different eigenvalues of the coefficients' matrix $\Theta$ denoted as $\lambda_1$ and $\lambda_2$, where $|\lambda_1|<1$, $|\lambda_2|<1$, and $1<\alpha<2$ we have 
			{\small%\begin{equation}
			\begin{align}
			{\rm{CV}}&(X_1(t),X_2(t-h))\nonumber\\&=\sum_{j=0}^{+\infty}\int_{S_2}\left(\frac{\lambda_2^ja_3s_1-\lambda_1^ja_3s_1+\lambda_2\lambda_1^js_2-\lambda_1\lambda_2^js_2+\lambda_2^ja_4s_2-\lambda_1^ja_4s_2}{\lambda_2-\lambda_1}\right)^{\langle \alpha-1 \rangle}\label{eq4}\\&\left(\frac{\lambda_1^h(\lambda_2\lambda_1^js_1-\lambda_1^ja_1s_1-\lambda_1^ja_2s_2)+\lambda_2^h(-\lambda_1\lambda_2^js_1+\lambda_2^ja_1s_1+\lambda_2^ja_2s_2)}{\lambda_2-\lambda_1}\right)\Gamma(ds).\nonumber
			\end{align}}
			%\end{equation}
			{\small%\begin{equation}
			\begin{align}
			{\rm{CV}}&(X_1(t),X_2(t+h))\nonumber\\&=\sum_{j=0}^{+\infty}\int_{S_2}\left(\frac{\lambda_2\lambda_1^js_1-\lambda_1\lambda_2^js_1+\lambda_2^ja_1s_1-\lambda_1^ja_1s_1+\lambda_2^ja_2s_2-\lambda_1^ja_2s_2}{\lambda_2-\lambda_1}\right)\label{eq8}\\&\left({\frac{\lambda_1^h(-\lambda_1^ja_3s_1+\lambda_2\lambda_1^js_2-\lambda_1^ja_4s_2)+\lambda_2^h(\lambda_2^ja_3s_1-\lambda_1\lambda_2^js_2+\lambda_2^ja_4s_2)}{\lambda_2-\lambda_1}}\right)^{\langle \alpha-1 \rangle}\Gamma(ds)\nonumber
			\end{align}}	
		%	\end{equation}
			\item in the case of equal eigenvalues of the coefficients' matrix $\Theta$ denoted as $\lambda_1=\lambda_2=\lambda$, where $|\lambda|<1$, and $1<\alpha<2$ we have 
			{\small%\begin{equation}
			\begin{align}
			{\rm{CV}}(X_1(t),X_2(t-h))=&\sum_{j=0}^{+\infty}\int_{S_2}\left(j\lambda^{j-1}a_3s_1-(j-1)\lambda^js_2+j\lambda^{j-1}a_4s_2\right)^{\langle \alpha-1 \rangle}\nonumber\\&\Big(\lambda^h{(j\lambda^{j-1}a_1s_1-(j-1)\lambda^js_1+j\lambda^{j-1}a_2s_2)}\label{eq6}\\&+h\lambda^h{(\lambda^{j-1}a_1s_1-\lambda^js_1+\lambda^{j-1}a_2s_2)}\Big)\Gamma(ds)\nonumber\\
			{\rm{CV}}(X_1(t),X_2(t+h))=&\sum_{j=0}^{+\infty}\int_{S_2}\left(j\lambda^{j-1}a_1s_1-(j-1)\lambda^{j}s_1+j\lambda^{j-1}a_2s_2\right)\nonumber\\&\Big(\lambda^h{(j\lambda^{j-1}a_3s_1-(j-1)\lambda^js_2+j\lambda^{j-1}a_4s_2)}\label{eq10}\\&+h\lambda^h{(\lambda^{j-1}a_3s_1+\lambda^{j-1}a_4s_2-\lambda^{j}s_2)}\Big)^{\langle \alpha-1 \rangle}\Gamma(ds)\nonumber
			\end{align}}
		%	\end{equation}}
		\end{enumerate}
\end{lemma}}
\begin{proof}
	The above equations follow from the formulas presented in the authors' previous paper, see \cite{nasza2}, and the expressions for the $j$-th power of the coefficients' matrix $\Theta$ given in Eq. (\ref{theta1}) or Eq. (\ref{theta2}).
\end{proof}
	\section{The asymptotic behavior for the ratio of codifference and covariation for bidimensional AR(1) model with $\alpha-$stable noise} \label{section_proof}
{In this section, we examine the asymptotic relation of the cross-codifference function and the cross-covariation function for the bidimensional AR(1) model with $\alpha-$stable noise presented in Section \ref{sec2}. Before formulating the relevant theorem, we separately consider the asymptotic behavior of both measures for $h \to +\infty$ in Lemma \ref{lema1} and Lemma \ref{lema2}, respectively. For both functions, we distinguish five separate cases, see Table 1.}
\begin{table}[htbp] \label{TABLE1} 
	\centering 
	{\normalsize{\caption{The considered cases for the asymptotic behavior of the cross-dependence measures.}}}\vspace{0.1cm}
	{\small \begin{tabular}{|m{1.6cm}|m{1.6cm}|m{1.9cm}|m{1.7cm}|m{1.7cm}|} \hline
			$|\lambda_1|>|\lambda_2|$ & $|\lambda_1|<|\lambda_2|$  & $\lambda_1=\lambda_2=\lambda$ & $\lambda_1=-\lambda_2$ and $h$ even & $\lambda_1=-\lambda_2$ and $h$ odd\\
			\hline\hline 
			CASE I & CASE II & CASE III & CASE IV & CASE V \\
			\hline
	\end{tabular}}
\end{table}

{\begin{lemma} \label{lema1}
		If $\{\mathbf{X}(t)\}=\{X_1(t),X_2(t)\}$ for $t\in Z$ is the bounded solution of Eq. (3) given by Eq. (6) with $1<\alpha<2$, then for the cross-codifference function the following asymptotic formulas are true when  $h \to +\infty$
		\begin{enumerate}[a)] 
			\item
			{\small\begin{align}
			{\rm{CD}}(X_1(t),X_2(t-h))\sim
			\left\{\begin{array}{ll}
			\alpha D_1\,\lambda_1^h \qquad \text{\rm{CASE I}},\\
			\alpha D_2\,\lambda_2^h \qquad \text{\rm{CASE II}},\\
			\alpha D_3\,h\lambda^h \qquad \text{\rm{CASE III}},\\
			\alpha (D_1+D_2)\,\lambda_1^h \qquad \text{\rm{CASE IV}},\\
			\alpha (D_1-D_2)\,\lambda_1^h \qquad \text{\rm{CASE V}},\nonumber
			\end{array} \right.
			\end{align}}
			where the constants $D_1$, $D_2$ and $D_3$ are given in Eqs. (\ref{D_1}), (\ref{D_2}) and (\ref{D_3}), respectively.
			\normalsize\item
			{\small\begin{align}
			{\rm{CD}}(X_1(t),X_2(t+h))\sim
			\left\{\begin{array}{ll}
			\alpha D_4\,\lambda_1^h \qquad \text{\rm{CASE I}},\\
			\alpha D_5\,\lambda_2^h \qquad \text{\rm{CASE II}},\\
			\alpha D_6\,h\lambda^h \qquad \text{\rm{CASE III}},\\
			\alpha (D_4+D_5)\,\lambda_1^h \qquad \text{\rm{CASE IV}},\\
			\alpha (D_4-D_5)\,\lambda_1^h \qquad \text{\rm{CASE V}},\nonumber
			\end{array} \right.
			\end{align}}
			where the constants $D_4$, $D_5$ and $D_6$ are given in Eqs. (\ref{D_4}), (\ref{D_5}) and (\ref{D_6}), respectively.
		\end{enumerate}
\end{lemma}}
{\begin{proof}
		The proof of Lemma \ref{lema1} is given in Appendix A.
\end{proof}}
{\begin{lemma} \label{lema2}
		If $\{\mathbf{X}(t)\}=\{X_1(t),X_2(t)\}$ for $t\in Z$ is the bounded solution of Eq. (3) given by Eq. (6) with $1<\alpha<2$, then for the cross-covariation function we have
		\begin{enumerate}[a)] 
			\item the following asymptotic formulas when $h \to +\infty$
			{\small\begin{align}
			{\rm{CV}}(X_1(t),X_2(t-h))\sim
			\left\{\begin{array}{ll}
			D_1\,\lambda_1^h \qquad \text{\rm{CASE I}},\\
			D_2\,\lambda_2^h \qquad \text{\rm{CASE II}},\\
			D_3\,h\lambda^h \qquad \text{\rm{CASE III}},\nonumber
			\end{array} \right.
			\end{align}}
			or the following exact formulas
			{\small\begin{align}
			{\rm{CV}}(X_1(t),X_2(t-h))=
			\left\{\begin{array}{ll}
			(D_1+D_2)\,\lambda_1^h \qquad \text{\rm{CASE IV}},\\
			(D_1-D_2)\,\lambda_1^h \qquad \text{\rm{CASE V}},\nonumber
			\end{array} \right.
			\end{align}}
			where the constants $D_1$, $D_2$ and $D_3$ are given in Eqs. (\ref{D_1}), (\ref{D_2}) and (\ref{D_3}), respectively.
			\normalsize\item the following asymptotic formulas when $h \to +\infty$
			{\small\begin{align}
			{\rm{CV}}(X_1(t),X_2(t+h))\sim
			\left\{\begin{array}{ll}
			D_7\,\left(\lambda_1^h\right)^{\langle \alpha-1 \rangle} \qquad \text{\rm{CASE I}},\\
			D_8\,\left(\lambda_2^h\right)^{\langle \alpha-1 \rangle} \qquad \text{\rm{CASE II}},\\
			D_9\,\left(h\lambda^h\right)^{\langle \alpha-1 \rangle} \qquad \text{\rm{CASE III}}\nonumber
			\end{array} \right.
			\end{align}}
			or the following exact formulas
			{\small\begin{align}
			{\rm{CV}}(X_1(t),X_2(t+h))=
			\left\{\begin{array}{ll}
			D_{10}\,\left(\lambda_1^h\right)^{\langle \alpha-1 \rangle} \qquad \text{\rm{CASE IV}},\\
			D_{11}\,\left(\lambda_1^h\right)^{\langle \alpha-1 \rangle} \qquad \text{\rm{CASE V}},\nonumber
			\end{array} \right.
			\end{align}}
			where the constants $D_7$, $D_8$, $D_9$, $D_{10}$ and $D_{11}$ are given in Eqs. (\ref{D_7}), (\ref{D_8}), (\ref{D_9}), (\ref{D_10}) and (\ref{D_11}), respectively.
		\end{enumerate}
\end{lemma}}
{\begin{proof} 
		The proof of Lemma \ref{lema2} is given in Appendix B.
\end{proof}}
{\begin{theorem} \label{theorem}
		If $\{\mathbf{X}(t)\}=\{X_1(t),X_2(t)\}$ is the bounded solution of Eq. (\ref{eq1}) given by Eq. (\ref{eq2}), then for $1<\alpha<2$ and $t\in Z$ the following formulas hold
		\begin{enumerate}[a)]	
			\begin{multicols}{2}
				\item
				{\small\begin{align}
				\lim_{h\to+\infty}\frac{{\rm{CD}}(X_1(t),X_2(t-h))}{{\rm{CV}}(X_1(t),X_2(t-h))}=\alpha,
				\end{align}}
			
				{assuming that ${\rm{CV}}(X_1(t),X_2(t-h)) \neq 0$.}
				\columnbreak
				\item
				{\small\begin{align}
				\lim_{h\to+\infty}\frac{{\rm{CD}}(X_1(t),X_2(t+h))}{{\rm{CV}}(X_1(t),X_2(t+h))}=0,
				\end{align}}
			
				{assuming that ${\rm{CV}}(X_1(t),X_2(t+h)) \neq 0$.}
			\end{multicols}
		\end{enumerate} 
		
\end{theorem}}
{\begin{proof}
		The limits given in Theorem \ref{theorem} follow directly from the formulas describing the asymptotic behavior of the cross-dependence measures given in Lemma \ref{lema1}-\ref{lema2}.
\end{proof}}

	\section{Example} \label{example_sec}
To illustrate the theoretical results showing the asymptotic behavior of the dependence measures for the bidimensional autoregressive model of order $1$ with general symmetric $\alpha-$stable noise, we consider two exemplary bidimensional AR(1) time series  with the following coefficients' matrices  
{\small\begin{equation*}
\Theta_1= \left[
\begin{array}{cc}
-0.2 & 0.1\\
-0.3 & 0.6
\end{array}
\right] \hspace{1cm} \text{and} \hspace{1cm} \Theta_2= \left[
\begin{array}{cc}
0.5 & 0.1\\
-0.1 & 0.7
\end{array}
\right].
\label{theta_example}
\end{equation*}}
Since the eigenvalues of the matrix $\Theta_1$ ($\lambda_1 \approx -0.16$, $\lambda_2 \approx 0.56$) and {$\Theta_2$ ($\lambda_1=\lambda_2=\lambda = 0.6$)} are less than $1$ in absolute value, for both time series there exists the bounded solution given in Eq. (\ref{eq2}). Moreover, let us assume that the spectral measure of the symmetric $\alpha-$stable random vector is concentrated on four points on the unit sphere $S_2$, namely it has the following form
{\small\begin{equation*}
\Gamma\left(\cdot\right)=0.5\,\delta\left(\left(\frac{1}{2},\frac{\sqrt{3}}{2}\right)\right)+0.5\,\delta\left(\left(-\frac{1}{2},-\frac{\sqrt{3}}{2}\right)\right)+0.2\,\delta\left(\left(-\frac{1}{2},\frac{\sqrt{3}}{2}\right)\right)+0.2\,\delta\left(\left(\frac{1}{2},-\frac{\sqrt{3}}{2}\right)\right).
\label{noise}
\end{equation*}}
\begin{figure}[h!] 
	\centering
	\includegraphics[scale=0.52]{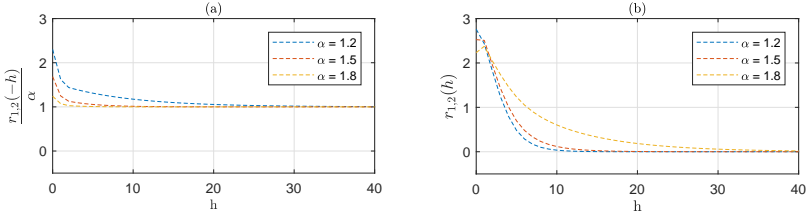}
	\caption{{The terms $\frac{r_{1,2}(-h)}{\alpha}$ (panel (a)) and ${r_{1,2}(h)}$ (panel (b)) presenting the asymptotic behavior of the cross-dependence measures for the bidimensional AR(1) model with coefficients given in matrix $\Theta_1$ and the $\alpha-$stable noise with the symmetric spectral measure specified in Eq. (\ref{noise}) for $\alpha=1.2$ (blue line), $\alpha=1.5$ (red line) and $\alpha=1.8$ (yellow line)}}
	\label{fig1}
\end{figure}
\begin{figure}[h!] 
	\centering
	\includegraphics[scale=0.52]{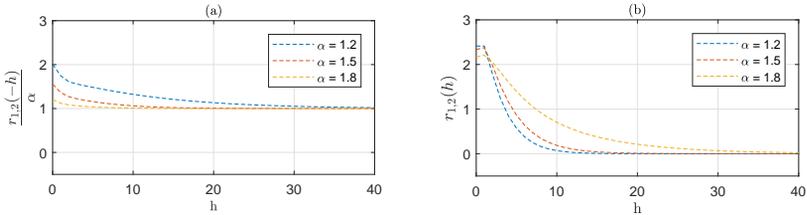}
	\caption{{The terms $\frac{r_{1,2}(-h)}{\alpha}$ (panel (a)) and ${r_{1,2}(h)}$ (panel (b)) presenting the asymptotic behavior of the dependence measures for the bidimensional AR(1) model with coefficients given in matrix $\Theta_2$ and the $\alpha-$stable noise with the symmetric spectral measure specified in Eq. (\ref{noise}) for $\alpha=1.2$ (blue line), $\alpha=1.5$ (red line) and $\alpha=1.8$ (yellow line)}}
	\label{fig2}
\end{figure}

This specific model was considered in the authors' previous paper, see Example 3.3 in \cite{nasza2}, where the formulas for the corresponding cross-dependence measures are presented. Now, in order to demonstrate the asymptotic relation between {the cross-codifference and the cross-covariation} we plot the functions
{\small
	\[ \begin{array}{lll}%
	\frac{r_{1,2}(-h)}{\alpha}=\frac{\rm{CD}(X_1(t),X_2(t-h))}{\alpha\,\rm{CV}(X_1(t),X_2(t-h))} 
	\end{array}\]%
	and
	\[ \begin{array}{lll}%
	r_{1,2}(h)=\frac{\rm{CD}(X_1(t),X_2(t+h))}{\rm{CV}(X_1(t),X_2(t+h))}
	\end{array}\]}
for $h=0,1,\ldots,40$ and $\alpha=1.2$, $\alpha=1.5$ and $\alpha=1.8$. The graphs are presented in Figure \ref{fig1} (the model with coefficients given in matrix $\Theta_1$) and \mbox{in Figure \ref{fig2}} (the model with coefficients given \mbox{in matrix $\Theta_2$}). According to the results presented in Section \ref{section_proof}, {the term $\frac{r_{1,2}(-h)}{\alpha}$ tends to $1$ (panel (a) of Figures \ref{fig1}-\ref{fig2}) and the term $r_{1,2}(h)$ tends to $0$ (panel (b) of Figures \ref{fig1}-\ref{fig2})}. Moreover, we can see that the rate of the convergence depends on the value of the parameter $\alpha$. {The quotient denoted as $r_{1,2}(-h)$ converges faster for larger values of the stability index, on the contrary to the quotient denoted as $r_{1,2}(h)$} for which the larger value the parameter $\alpha$ takes, the slower is the convergence. 

\section{Conclusions} \label{Conclusions}

In this paper, we considered a bidimensional autoregressive model of order $1$ with $\alpha-$stable noise. For this process, the structure of dependence can be described using the codifference function and the covariation function. Here, we examined the asymptotic relation between these two measures {applied to the description of} the cross-dependence between the components of the bidimensional model (the cross-codifference function and the cross-covariation function). The main results concerning the asymptotic relationship between the considered dependence measures are presented in Theorem \ref{theorem} where we identified the case when the measures are asymptotically proportional with the coefficient of proportionality equal to the parameter $\alpha$. The conclusions obtained here constitute an extension of the results presented in \cite{arma} for the one-dimensional autoregressive models and may be useful in the context of the parameter $\alpha$ estimation. This can be a starting point for the introduction of a new estimation algorithm for the stability index for the bidimensional autoregressive model of order $1$ with the $\alpha-$stable noise.

	\appendix 
\section*{Appendix A: {Proof of Lemma \ref{lema1}}}
	{\begin{enumerate}[a)]
			\item First, we examine the asymptotic behavior of ${\rm{CD}}(X_1(t),X_2(t-h))$ for $h \to +\infty$. 
			\begin{itemize}
				\item Let us consider the case of two different real eigenvalues of the coefficients' matrix $\Theta$, $\lambda_1 \neq \lambda_2$, and $|\lambda_1|<1$, $|\lambda_2|<1$.
				For the cross-codifference function ${\rm{CD}}(X_1(t),X_2(t-h))$ given in Lemma \ref{th1} (see Eq. (\ref{eq3})) we introduce the notation
				{\footnotesize %\begin{equation}
					\begin{align}
					A_1(s_1,s_2,a_1,a_2,\lambda_1,\lambda_2,j)&=A_1=\frac{\lambda_2\lambda_1^js_1-\lambda_1^ja_1s_1-\lambda_1^ja_2s_2}{\lambda_2-\lambda_1},\nonumber\\
					B_1(s_1,s_2,a_1,a_2,\lambda_1,\lambda_2,j)&=B_1=\frac{-\lambda_1\lambda_2^js_1+\lambda_2^ja_1s_1+\lambda_2^ja_2s_2}{\lambda_2-\lambda_1},\label{A1B1C1}\\
					C_1(s_1,s_2,a_3,a_4,\lambda_1,\lambda_2,j)&=C_1=\frac{\lambda_1^j(-a_3s_1+\lambda_2s_2-a_4s_2)+\lambda_2^j(a_3s_1-\lambda_1s_2+a_4s_2)}{\lambda_2-\lambda_1}.\nonumber
					\end{align}}%\end{equation}}
				In the following part of the proof, the above expressions will be referred as $A_1$, $B_1$ and $C_1$ to simplify the notation. However, it is important to notice that although they do not depend on h, they are dependent on $j$, $s_1$, $s_2$, the coefficients $a_1,a_2,a_3,a_4$ and thus by the eigenvalues of matrix $\Theta$.
				Now, the cross-codifference function ${\rm{CD}}(X_1(t),X_2(t-h))$ given in Eq. (\ref{eq3}) takes the form
				{\footnotesize \begin{multline}
					{\rm{CD}}(X_1(t),X_2(t-h))\\={\sum_{j=0}^{+\infty}\int_{S_2}\left(|\lambda_1^h{A_1}+\lambda_2^h{B}|^{\alpha}+|C_1|^{\alpha}-|C_1-(\lambda_1^h{A_1}+\lambda_2^h{B_1})|^\alpha\right)\Gamma(ds)}.
					\label{p1}
					\end{multline}}
				\begin{enumerate}[I)]
					\item Let us consider the case of $|\lambda_1|>|\lambda_2|$. Now, we split the proof into two parts. Namely, we show that
					\begin{enumerate}[i)]
						\item 
						{\footnotesize \begin{multline} \label{equality1}
							\hspace{-0.5cm}\lim_{h \to +\infty}\sum_{j=0}^{+\infty}\int_{S_2}\frac{\left|\lambda_1^hA_1+\lambda_2^hB_1\right|^{\alpha}+\left|C_1\right|^{\alpha}-\left|C_1-(\lambda_1^hA_1+\lambda_2^hB_1)\right|^\alpha}{\lambda_1^h}\Gamma(ds)\\\hspace{-0.5cm}=\sum_{j=0}^{+\infty}\lim_{h \to +\infty}\int_{S_2}\frac{\left|\lambda_1^hA_1+\lambda_2^hB_1\right|^{\alpha}+\left|C_1\right|^{\alpha}-\left|C_1-(\lambda_1^hA_1+\lambda_2^hB_1)\right|^\alpha}{\lambda_1^h}\Gamma(ds),
							\end{multline}}
						\item
						{\footnotesize \begin{multline} \label{equality2}
							\hspace{-0.5cm}\lim_{h \to +\infty}\int_{S_2}\frac{\left|\lambda_1^hA_1+\lambda_2^hB_1\right|^{\alpha}+\left|C_1\right|^{\alpha}-\left|C_1-(\lambda_1^hA_1+\lambda_2^hB_1)\right|^\alpha}{\lambda_1^h}\Gamma(ds)\\\hspace{-0.5cm}=\int_{S_2}\lim_{h \to +\infty}\frac{\left|\lambda_1^hA_1+\lambda_2^hB_1\right|^{\alpha}+\left|C_1\right|^{\alpha}-\left|C_1-(\lambda_1^hA_1+\lambda_2^hB_1)\right|^\alpha}{\lambda_1^h}\Gamma(ds).
							\end{multline}}
					\end{enumerate}
					Let us focus on part i). According to the dominated convergence theorem \cite[Ch.~5]{dom_con_theorem}, the equality given in Eq. (\ref{equality1}) is true if the infinite sum over $j$ converges uniformly. Let us notice that using the following inequalities true for all $ n,m \in \R, 1<\alpha<2$ \cite{maejima2003}
					{\footnotesize \begin{equation}
						\begin{aligned} \label{formula1}
						\left||n|^\alpha+|m|^\alpha-|n+m|^\alpha\right|&\leq(\alpha+1)|n|^\alpha+\alpha|n||m|^{\alpha-1}\\
						\left|n+m\right|^\alpha&\leq2^{\alpha-1}\left(|n|^\alpha+|m|^\alpha\right),
						\end{aligned}
						\end{equation}}
					one can show that for all $j \in \N$ we have
					{\footnotesize \begin{align*}
						&\left|\int_{S_2}\frac{\left|\lambda_1^hA_1+\lambda_2^hB_1\right|^{\alpha}+\left|C_1\right|^{\alpha}-\left|C_1-(\lambda_1^hA_1+\lambda_2^hB_1)\right|^\alpha}{\lambda_1^h}\Gamma(ds)\right|\\ \leq&\int_{S_2}\left|\frac{\left|\lambda_1^hA_1+\lambda_2^hB_1\right|^{\alpha}+\left|C_1\right|^{\alpha}-\left|C_1-(\lambda_1^hA_1+\lambda_2^hB_1)\right|^\alpha}{\lambda_1^h}\right|\Gamma(ds)\\\leq& \int_{S_2}\frac{(\alpha+1)\left|\lambda_1^hA_1+\lambda_2^hB_1\right|^{\alpha}+\alpha\left|\lambda_1^hA_1+\lambda_2^hB_1\right||C_1|^{\alpha-1}}{|\lambda_1^h|}\Gamma(ds) \\ = &\,(\alpha+1)\left|\lambda_1^h\right|^{\alpha-1}\int_{S_2}\left|A_1+\left(\frac{\lambda_2}{\lambda_1}\right)^hB_1\right|^{\alpha}\Gamma(ds)\\+&\,\alpha\int_{S_2}\left|A_1+\left(\frac{\lambda_2}{\lambda_1}\right)^hB_1\right||C_1|^{\alpha-1}\Gamma(ds) \\ \leq & \,(\alpha+1)\left|\lambda_1^h\right|^{\alpha-1}2^{\alpha-1}\int_{S_2}\left(\left|A_1\right|^\alpha+\left|\left(\frac{\lambda_2}{\lambda_1}\right)^hB_1\right|^{\alpha}\right)\Gamma(ds)\\+&\,\alpha\int_{S_2}\left(\left|A_1\right|+\left|\left(\frac{\lambda_2}{\lambda_1}\right)^hB_1\right|\right)|C_1|^{\alpha-1}\Gamma(ds)\\ \leq & \,2^{\alpha-1}(\alpha+1)\left(\int_{S_2}\left|A_1\right|^\alpha\Gamma(ds)+\int_{S_2}\left|B_1\right|^{\alpha}\Gamma(ds)\right)\\+&\alpha\left(\int_{S_2}\left|A_1\right||C_1|^{\alpha-1}\Gamma(ds)+\int_{S_2}\left|B_1\right||C_1|^{\alpha-1}\Gamma(ds)\right)=M_j,
						\end{align*}}
					and $M_j$ is independent of $h$. Now, the infinite sum under the limit in Eq. (\ref{equality1}) converges uniformly if the infinite sum of $M_j$ over $j \in \N$ converges, which is equivalent to the following set of conditions
					{\footnotesize 
						\begin{equation}
						\begin{aligned}\label{conditions1} 
						&\sum_{j=0}^{+\infty}\int_{S_2}\left|A_1\right|^\alpha\Gamma(ds)<+\infty, \quad \sum_{j=0}^{+\infty}\int_{S_2}\left|A_1\right||C_1|^{\alpha-1}\Gamma(ds)<+\infty,\\
						&\sum_{j=0}^{+\infty}\int_{S_2}\left|B_1\right|^{\alpha}\Gamma(ds)<+\infty, \quad
						\sum_{j=0}^{+\infty}\int_{S_2}\left|B_1\right||C_1|^{\alpha-1}\Gamma(ds)<+\infty,
						\end{aligned}
						\end{equation}}
					{which are always satisfied (see Remark \ref{remark} in Appendix C) and consequently the equality given in Eq. (\ref{equality1}) is true.} Let us focus on part ii). Similarly,  according to the dominated convergence theorem, the equality given in Eq. (\ref{equality2}) is satisfied if the integrand is dominated by an integrable function independent of $h$. Let us notice that using the formulas given in Eq. (\ref{formula1}) for a fixed $\mathbf{s}=(s_1,s_2)\in S_2$ we have that
					{\footnotesize \begin{align*}
						&\left|\frac{\left|\lambda_1^hA_1+\lambda_2^hB_1\right|^{\alpha}+\left|C_1\right|^{\alpha}-\left|C_1-(\lambda_1^hA_1+\lambda_2^hB_1)\right|^\alpha}{\lambda_1^h}\right|\\ \leq& \frac{(\alpha+1)\left|\lambda_1^hA_1+\lambda_2^hB_1\right|^{\alpha}+\alpha\left|\lambda_1^hA_1+\lambda_2^hB_1\right||C_1|^{\alpha-1}}{|\lambda_1^h|} \\ = &\,(\alpha+1)\left|\lambda_1^h\right|^{\alpha-1}\left|A_1+\left(\frac{\lambda_2}{\lambda_1}\right)^hB_1\right|^{\alpha}+\,\alpha\left|A_1+\left(\frac{\lambda_2}{\lambda_1}\right)^hB_1\right||C_1|^{\alpha-1} \\ \leq & \,(\alpha+1)\left|\lambda_1^h\right|^{\alpha-1}2^{\alpha-1}\left(\left|A_1\right|^\alpha+\left|\left(\frac{\lambda_2}{\lambda_1}\right)^hB_1\right|^{\alpha}\right)\\+&\,\alpha\left(\left|A_1\right|+\left|\left(\frac{\lambda_2}{\lambda_1}\right)^hB_1\right|\right)|C_1|^{\alpha-1}\\ \leq & \,2^{\alpha-1}(\alpha+1)\left(\left|A_1\right|^\alpha+\left|B_1\right|^{\alpha}\right)+\alpha\left(\left|A_1\right||C_1|^{\alpha-1}+\left|B_1\right||C_1|^{\alpha-1}\right),
						\end{align*}}
					which is independent of $h$. Since the dominating function should be integrable we obtain that for all $j \in \N$ the following set of conditions has to be satisfied
					{\footnotesize 
						\begin{equation*}
						\begin{aligned} 
						&\int_{S_2}\left|A_1\right|^\alpha\Gamma(ds)<+\infty, \quad \int_{S_2}\left|A_1\right||C_1|^{\alpha-1}\Gamma(ds)<+\infty,\\
						&\int_{S_2}\left|B_1\right|^{\alpha}\Gamma(ds)<+\infty, \quad
						\int_{S_2}\left|B_1\right||C_1|^{\alpha-1}\Gamma(ds)<+\infty.
						\end{aligned}
						\end{equation*}}
					{Let us notice that the above conditions are satisfied since the conditions given in Eq. (\ref{conditions1}) hold.} Now, we calculate the limit of the integrand given in Eq. (\ref{equality2}) for a fixed $\mathbf{s}=(s_1,s_2)\in S_2$, namely
					{\footnotesize \begin{align*} 
						\lim_{h \to +\infty}\frac{\left|\lambda_1^hA_1+\lambda_2^hB_1\right|^{\alpha}+\left|C_1\right|^{\alpha}-\left|C_1-(\lambda_1^hA_1+\lambda_2^hB_1)\right|^\alpha}{\lambda_1^h}.
						\end{align*}}
					Let us notice that for $h \to +\infty$ we have
					{\footnotesize \begin{equation*}
						\begin{aligned}
						&\left|\lambda_1^h{A_1}+\lambda_2^h{B_1}\right|^{\alpha}+\left|C_1\right|^{\alpha}-\left|C_1-\left(\lambda_1^h{A_1}+\lambda_2^h{B_1}\right)\right|^\alpha\\=&\left|\lambda_1^h\left({A_1}+\left({\lambda_2}/{\lambda_1}\right)^h{B_1}\right)\right|^{\alpha}+\left|C_1\right|^{\alpha}-\left|C_1-\lambda_1^h\left({A_1}+\left(\lambda_2/\lambda_1\right)^h{B_1}\right)\right|^\alpha\\\sim&\left|\lambda_1^h\right|^{\alpha}\left|A_1\right|^{\alpha}+\left|C_1\right|^{\alpha}-\left|C_1-\lambda_1^hA_1\right|^\alpha
						\end{aligned}
						\end{equation*}}
					and since for $1<\alpha<2$ and $a,c \in \R$ the following limit holds 
					{\footnotesize \begin{align}
						\lim_{x \to 0}\frac{|ax|^\alpha+|c|^\alpha-|c-ax|^\alpha}{x}=\alpha\, ac^{\langle\alpha-1\rangle}
						\label{limit}
						\end{align}} 
					and $\lambda_1^h \to 0$ \rm{for} $h \to +\infty$, we obtain
					{\footnotesize \begin{align*}
						\lim_{h \to +\infty}\frac{\left|\lambda_1^hA_1\right|^{\alpha}+\left|C_1\right|^{\alpha}-\left|C_1-\lambda_1^hA_1\right|^\alpha}{\lambda_1^h}=\alpha\, A_1C_1^{\langle\alpha-1\rangle}
						\end{align*}}
					for a fixed $\mathbf{s}=(s_1,s_2)\in S_2$. {Now since the equalities presented in i) and ii) are true, we obtain}
					{\footnotesize \begin{multline}
						\lim_{h \to +\infty}\sum_{j=0}^{+\infty}\int_{S_2}\frac{\left|\lambda_1^hA_1+\lambda_2^hB_1\right|^{\alpha}+\left|C_1\right|^{\alpha}-\left|C_1-(\lambda_1^hA_1+\lambda_2^hB_1)\right|^\alpha}{\lambda_1^h}\Gamma(ds)\\=\alpha D_1,\label{p2}
						\end{multline}}
					which is equivalent to the fact that
					{\footnotesize \begin{multline}
						\sum_{j=0}^{+\infty}\int_{S_2}{\left|\lambda_1^hA_1+\lambda_2^hB_1\right|^{\alpha}+\left|C_1\right|^{\alpha}-\left|C_1-(\lambda_1^hA_1+\lambda_2^hB_1)\right|^\alpha}\Gamma(ds)\\\sim\alpha D_1\,\lambda_1^h\ \ \text{\normalsize for}\ \ h \to +\infty,
						\label{p3}
						\end{multline}}
					where
					{\footnotesize \begin{align}\label{D_1}
						D_1:=\sum_{j=0}^{+\infty}\int_{S_2} A_1C_1^{\langle\alpha-1\rangle}\Gamma(ds).
						\end{align}}
					Let us notice that the conditions given in Eq. (\ref{conditions1}) guarantee that $D_1<+\infty$.
					\item Let us consider the case when $|\lambda_1|<|\lambda_2|$. {Similary as above we have}
					{\footnotesize \begin{multline}
						\sum_{j=0}^{+\infty}\int_{S_2}{\left|\lambda_1^hA_1+\lambda_2^hB_1\right|^{\alpha}+\left|C_1\right|^{\alpha}-\left|C_1-(\lambda_1^hA_1+\lambda_2^hB_1)\right|^\alpha}\Gamma(ds)\\\sim\alpha D_2\,\lambda_2^h\ \ \text{\normalsize for}\ \ h \to +\infty,
						\end{multline}}
					where
					{\footnotesize \begin{align}\label{D_2}
						D_2:=\sum_{j=0}^{+\infty}\int_{S_2}B_1C_1^{\langle\alpha-1\rangle}\Gamma(ds)<+\infty.
						\end{align}}
					\item Let us consider the case when $\lambda_1=-\lambda_2$ and $h$ is even. Proceeding in the same way as above, we obtain that
					{\footnotesize \begin{multline}
						\sum_{j=0}^{+\infty}\int_{S_2}{\left|\lambda_1^hA_1+\lambda_2^hB_1\right|^{\alpha}+\left|C_1\right|^{\alpha}-\left|C_1-(\lambda_1^hA_1+\lambda_2^hB_1)\right|^\alpha}\Gamma(ds)\\\sim\alpha (D_1+D_2)\,\lambda_1^h\ \ \text{\normalsize for}\ h \to +\infty,
						\end{multline}}
					where
					{\footnotesize \begin{align*}
						D_1+D_2=\sum_{j=0}^{+\infty}\int_{S_2}(A_1+B_1)C_1^{\langle\alpha-1\rangle}\Gamma(ds)<+\infty. \end{align*}}
					\item Let us consider the case when $\lambda_1=-\lambda_2$ and $h$ is odd. Similarly to the previous cases, we have
					{\footnotesize \begin{multline}
						\sum_{j=0}^{+\infty}\int_{S_2}{\left|\lambda_1^hA_1+\lambda_2^hB_1\right|^{\alpha}+\left|C_1\right|^{\alpha}-\left|C_1-(\lambda_1^hA_1+\lambda_2^hB_1)\right|^\alpha}\Gamma(ds)\\\sim\alpha (D_1-D_2)\,\lambda_1^h\ \text{\normalsize for}\ h \to +\infty,
						\end{multline}}
					where
					{\footnotesize \begin{align*}
						D_1-D_2=\sum_{j=0}^{+\infty}\int_{S_2}(A_1-B_1)C_1^{\langle\alpha-1\rangle}\Gamma(ds)<+\infty.
						\end{align*}}
				\end{enumerate}
				\item Let us consider the case of the real and equal eigenvalues of the coefficients' matrix $\Theta$, $\lambda_1=\lambda_2=\lambda$, and $|\lambda|<1$. For the codifference $\rm{CD}(X_1(t),X_2(t-h))$ given in Lemma \ref{th1} (see Eq. (\ref{eq5})) we introduce the notation
				{\footnotesize \begin{equation}\begin{aligned} \label{A2B2C2}
					A_2(s_1,s_2,a_1,a_2,\lambda,j)&=A_2=j\lambda^{j-1}a_1s_1-j\lambda^js_1+\lambda^js_1+j\lambda^{j-1}a_2s_2,\\
					B_2(s_1,s_2,a_1,a_2,\lambda,j)&=B_2=\lambda^{j-1}a_1s_2-\lambda^js_1+\lambda^{j-1}a_2s_2,\\
					C_2(s_1,s_2,a_3,a_4,\lambda,j)&=C_2=j\lambda^{j-1}a_3s_1+j\lambda^{j-1}a_4s_2-j\lambda^js_2+\lambda^js_2.	
					\end{aligned}\end{equation}}
				Again, the above expressions will be referred to as $A_2$, $B_2$ and $C_2$ to simplify the notation. Now, the cross-codifference function $\rm{CD}(X_1(t),X_2(t-h))$ given in Eq. (\ref{eq5}) takes the form
				{\footnotesize \begin{multline}
					{\rm{CD}}(X_1(t),X_2(t-h))\\=\sum_{j=0}^{+\infty}\int_{S_2}\left(|	\lambda^h{A_2}+h\lambda^h{B_2}|^{\alpha}+|C_2|^{\alpha}-|C_2-(\lambda^h{A_2}+h\lambda^h{B_2})|^\alpha\right)\Gamma(ds).
					\label{p4}
					\end{multline}}
				Now, similarly to the case of two different eigenvalues, we split the proof into two parts. Namely, we want to show that
				\begin{enumerate}[i)]
					\item 
					{\footnotesize \begin{multline} \label{equality3}
						\hspace{-0.5cm}\lim_{h \to +\infty}\sum_{j=0}^{+\infty}\int_{S_2}\frac{\left|\lambda^hA_2+h\lambda^hB_2\right|^{\alpha}+\left|C_2\right|^{\alpha}-\left|C_2-(\lambda^hA_2+h\lambda^hB_2)\right|^\alpha}{h\lambda^h}\Gamma(ds)\\\hspace{-0.5cm}=\sum_{j=0}^{+\infty}\lim_{h \to +\infty}\int_{S_2}\frac{\left|\lambda^hA_2+h\lambda^hB_2\right|^{\alpha}+\left|C_2\right|^{\alpha}-\left|C_2-(\lambda^hA_2+h\lambda^hB_2)\right|^\alpha}{h\lambda^h}\Gamma(ds),
						\end{multline}}
					\item
					{\footnotesize \begin{multline} \label{equality4}
						\hspace{-0.5cm}\lim_{h \to +\infty}\int_{S_2}\frac{\left|\lambda^hA_2+h\lambda^hB_2\right|^{\alpha}+\left|C_2\right|^{\alpha}-\left|C_2-(\lambda^hA_2+h\lambda^hB_2)\right|^\alpha}{h\lambda^h}\Gamma(ds)\\\hspace{-0.5cm}=\int_{S_2}\lim_{h \to +\infty}\frac{\left|\lambda^hA_2+h\lambda^hB_2\right|^{\alpha}+\left|C_2\right|^{\alpha}-\left|C_2-(\lambda^hA_2+h\lambda^hB_2)\right|^\alpha}{h\lambda^h}\Gamma(ds).
						\end{multline}}
				\end{enumerate}
				To justify the equality in part i), according to the dominated convergence theorem, the infinite sum over $j$ has to converge uniformly. Using the inequalities given in Eq. (\ref{formula1}) we can show that for all $j \in \N$ we have
				{\footnotesize \begin{align*}
					&\left|\int_{S_2}\frac{\left|\lambda^hA_2+h\lambda^hB_2\right|^{\alpha}+\left|C_2\right|^{\alpha}-\left|C_2-(\lambda^hA_2+h\lambda^hB_2)\right|^\alpha}{h\lambda^h}\Gamma(ds)\right| \\ \leq&\int_{S_2}\left|\frac{\left|\lambda^hA_2+h\lambda^hB_2\right|^{\alpha}+\left|C_2\right|^{\alpha}-\left|C_2-(\lambda^hA_2+h\lambda^hB_2)\right|^\alpha}{h\lambda^h}\right|\Gamma(ds) \\ \leq &\int_{S_2}\frac{(\alpha+1)\left|\lambda^hA_2+h\lambda^hB_2\right|^{\alpha}+\alpha\left|\lambda^hA_2+h\lambda^hB_2\right||C_2|^{\alpha-1}}{|h\lambda^h|}\Gamma(ds) \\ = &\,(\alpha+1)\left|h\lambda^h\right|^{\alpha-1}\int_{S_2}\left|\frac{A_2}{h}+B_2\right|^{\alpha}\Gamma(ds)+\,\alpha\int_{S_2}\left|\frac{A_2}{h}+B_2\right||C_2|^{\alpha-1}\Gamma(ds) \\ \leq & \,(\alpha+1)\left|h\lambda^h\right|^{\alpha-1}2^{\alpha-1}\int_{S_2}\left(\left|\frac{A_2}{h}\right|^\alpha+\left|B_2\right|^{\alpha}\right)\Gamma(ds)\\+&\,\alpha\int_{S_2}\left(\left|\frac{A_2}{h}\right|+\left|B_2\right|\right)|C|^{\alpha-1}\Gamma(ds)\\ \leq & \,2^{\alpha-1}(\alpha+1)M\left(\int_{S_2}\left|A_2\right|^\alpha\Gamma(ds)+\int_{S_2}\left|B_2\right|^{\alpha}\Gamma(ds)\right)\\+&\alpha\left(\int_{S_2}\left|A_2\right||C_2|^{\alpha-1}\Gamma(ds)+\int_{S_2}\left|B_2\right||C_2|^{\alpha-1}\Gamma(ds)\right)=N_j,
					\end{align*}}
				where $M$ denotes the boundary of the sequence $\{h\lambda^h\}$ for $h\in\N$ (the sequence converges). Let us notice that $N_j$ is independent of $h$. Now, the infinite sum under the limit in Eq. (\ref{equality3}) converges uniformly if the infinite sum of $N_j$ over $j \in \N$ converges, which is equivalent to the following set of conditions
				{\footnotesize \begin{equation}
					\begin{aligned}\label{conditions3} 
					&\sum_{j=0}^{+\infty}\int_{S_2}\left|A_2\right|^\alpha\Gamma(ds)<+\infty, \quad \sum_{j=0}^{+\infty}\int_{S_2}\left|A_2\right||C_2|^{\alpha-1}\Gamma(ds)<+\infty,\\
					&\sum_{j=0}^{+\infty}\int_{S_2}\left|B_2\right|^{\alpha}\Gamma(ds)<+\infty, \quad
					\sum_{j=0}^{+\infty}\int_{S_2}\left|B_2\right||C_2|^{\alpha-1}\Gamma(ds)<+\infty.
					\end{aligned}\end{equation}}	
				{The above conditions are always satisfied (see Remark \ref{remark} in Appendix C) and consequently the equality given in Eq. (\ref{equality3}) is true.} Let us focus on part ii). Similarly,  according to the dominated convergence theorem, the equality given in Eq. (\ref{equality4}) is satisfied if the integrand is dominated by an integrable function independent of $h$. Let us notice that using the inequalities given in Eq. (\ref{formula1}) for a fixed $\mathbf{s}=(s_1,s_2)\in S_2$ we have that
				{\footnotesize \begin{align*}
					&\left|\frac{\left|\lambda^hA_2+h\lambda^hB_2\right|^{\alpha}+\left|C_2\right|^{\alpha}-\left|C_2-(\lambda^hA_2+h\lambda^hB_2)\right|^\alpha}{h\lambda^h}\right| \\ \leq& \frac{(\alpha+1)\left|\lambda^hA_2+h\lambda^hB_2\right|^{\alpha}+\alpha\left|\lambda^hA_2+h\lambda^hB_2\right||C_2|^{\alpha-1}}{|h\lambda^h|} \\ = &\,(\alpha+1)\left|h\lambda^h\right|^{\alpha-1}\left|\frac{A_2}{h}+B_2\right|^{\alpha}+\,\alpha\left|\frac{A_2}{h}+B_2\right||C|^{\alpha-1} \\ \leq & \,(\alpha+1)\left|h\lambda^h\right|^{\alpha-1}2^{\alpha-1}\left(\left|\frac{A_2}{h}\right|^\alpha+\left|B_2\right|^{\alpha}\right)+\,\alpha\left(\left|\frac{A_2}{h}\right|+\left|B_2\right|\right)|C_2|^{\alpha-1}\\ \leq & \,2^{\alpha-1}(\alpha+1)M\left(\left|A_2\right|^\alpha+\left|B_2\right|^{\alpha}\right)+\alpha\left(\left|A_2\right||C_2|^{\alpha-1}+\left|B_2\right||C_2|^{\alpha-1}\right),
					\end{align*}}
				which is independent of $h$. Since the dominating function should be integrable we obtain that for all $j \in \N$ the following set of conditions has to be satisfied
				{\footnotesize \begin{equation*}
					\begin{aligned}
					&\int_{S_2}\left|A_2\right|^\alpha\Gamma(ds)<+\infty, \quad \int_{S_2}\left|A_2\right||C_2|^{\alpha-1}\Gamma(ds)<+\infty,\\
					&\int_{S_2}\left|B_2\right|^{\alpha}\Gamma(ds)<+\infty, \quad
					\int_{S_2}\left|B_2\right||C_2|^{\alpha-1}\Gamma(ds)<+\infty.
					\end{aligned}\end{equation*}}
				{Let us notice that the above conditions are satisfied since the conditions given in Eq. (\ref{conditions3}) hold.} Now, we calculate the limit of the integrand given in Eq. (\ref{equality4}) for a fixed $\mathbf{s}=(s_1,s_2)\in S_2$, namely
				{\footnotesize \begin{align*} 
					\lim_{h \to +\infty}\frac{\left|\lambda^hA_2+h\lambda^hB_2\right|^{\alpha}+\left|C_2\right|^{\alpha}-\left|C_2-(\lambda^hA_2+h\lambda^hB_2)\right|^\alpha}{h\lambda^h}.
					\end{align*}}
				Let us notice that for $h \to +\infty$ we have
				{\footnotesize \begin{equation*} \begin{aligned}
					&\left|\lambda^h{A_2}+h\lambda^h{B_2}\right|^{\alpha}+\left|C_2\right|^{\alpha}-\left|C_2-\left(\lambda^h{A_2}+h\lambda^h{B_2}\right)\right|^\alpha\\=&\left|h\lambda^h\left({A_2/h}+{B_2}\right)\right|^{\alpha}+\left|C_2\right|^{\alpha}-\left|C_2-h\lambda^h\left({A_2/h}+{B_2}\right)\right|^\alpha\\\sim&\left|h\lambda^h\right|^{\alpha}\left|B_2\right|^{\alpha}+\left|C_2\right|^{\alpha}-\left|C_2-h\lambda^hB_2\right|^\alpha.
					\end{aligned}\end{equation*}}
				Now, using the limit given in Eq. (\ref{limit}) and since $h\lambda^h \to 0$ \rm{for} $h \to +\infty$, we obtain
				{\footnotesize \begin{align*}
					\lim_{h \to +\infty}\frac{\left|h\lambda^hB_2\right|^{\alpha}+\left|C_2\right|^{\alpha}-\left|C_2-h\lambda^hB_2\right|^\alpha}{h\lambda^h}=\alpha\, B_2C_2^{\langle\alpha-1\rangle}
					\end{align*}}
				for a fixed $\mathbf{s}=(s_1,s_2)\in S_2$. {Now since the equalities presented in i) and ii) are true, we have}
				{\footnotesize \begin{multline}
					\lim_{h \to +\infty}\sum_{j=0}^{+\infty}\int_{S_2}\frac{\left|\lambda^hA_2+h\lambda^hB_2\right|^{\alpha}+\left|C_2\right|^{\alpha}-\left|C_2-(\lambda^hA_2+h\lambda^hB_2)\right|^\alpha}{h\lambda^h}\Gamma(ds)\\=\alpha D_3,
					\label{p6}
					\end{multline}}
				which is equivalent to the fact that
				{\footnotesize \begin{multline}
					\sum_{j=0}^{+\infty}\int_{S_2} \left|\lambda^hA_2+h\lambda^hB_2\right|^{\alpha}+\left|C_2\right|^{\alpha}-\left|C_2-(\lambda^hA_2+h\lambda^hB_2)\right|^\alpha\Gamma(ds)\\ \sim \alpha D_3\,h\lambda^h \ \ \text{\normalsize for}\ \ h \to +\infty,
					\label{p7}
					\end{multline}}
				where
				{\footnotesize{\begin{align} \label{D_3}
					D_3:=\sum_{j=0}^{+\infty}\int_{S_2} B_2C_2^{\langle\alpha-1\rangle}\Gamma(ds).
					\end{align}}}
				\normalsize Let us notice that the conditions given in Eq. (\ref{conditions3}) guarantee that $D_3<+\infty$.
			\end{itemize}
			\item Now, we examine the asymptotic behavior of ${\rm{CD}}(X_1(t),X_2(t+h))$ for $h \to +\infty$. 
			\begin{itemize}
				\item Let us consider the case of two different real eigenvalues of the coefficients' matrix $\Theta$, $\lambda_1 \neq \lambda_2$, and $|\lambda_1|<1$, $|\lambda_2|<1$. For the cross-codifference function ${\rm{CD}}(X_1(t),X_2(t+h))$ given in Lemma \ref{th1} (see Eq. (\ref{eq7})) we introduce the notation
				{\footnotesize %\begin{equation}
					\begin{align}
					A_3(s_1,s_2,a_3,a_4,\lambda_1,\lambda_2,j)&=A_3=\frac{-\lambda_1^ja_3s_1+\lambda_2\lambda_1^js_2-\lambda_1^ja_4s_2}{\lambda_2-\lambda_1},\nonumber\\
					B_3(s_1,s_2,a_3,a_4,\lambda_1,\lambda_2,j)&=B_3=\frac{\lambda_2^ja_3s_1-\lambda_1\lambda_2^js_2+\lambda_2^ja_4s_2}{\lambda_2-\lambda_1},\label{A3B3C3}\\
					C_3(s_1,s_2,a_1,a_2,\lambda_1,\lambda_2,j)&=C_3=\frac{\lambda_1^j(\lambda_2s_1-a_1s_1-a_2s_2)+\lambda_2^j(-\lambda_1s_1+a_1s_1+a_2s_2)}{\lambda_2-\lambda_1}.\nonumber
					\end{align}}%\end{equation}}
				Again, the above expressions will be referred to as $A_3$, $B_3$ and $C_3$ to simplify the notation. Now, we can write the formula for ${\rm{CD}}(X_1(t),X_2(t+h))$ given Eq. (\ref{eq7}) in the following form
				{\footnotesize  \begin{multline}
					{\rm{CD}}(X_1(t),X_2(t+h))\\=\sum_{j=0}^{+\infty}\int_{S_2}\left(|	\lambda_1^h{A_3}+\lambda_2^h{B_3}|^{\alpha}+|C_3|^{\alpha}-|C_3-(\lambda_1^h{A_3}+\lambda_2^h{B_3})|^\alpha\right)\Gamma(ds).
					\end{multline}}
				\begin{enumerate}[I)]
					\item Let us consider the case of $|\lambda_1|>|\lambda_2|$. Proceeding exactly as in a) we obtain
					{\footnotesize \begin{multline}
						\sum_{j=0}^{+\infty}\int_{S_2}{\left|\lambda_1^hA_3+\lambda_2^hB_3\right|^{\alpha}+\left|C_3\right|^{\alpha}-\left|C_3-(\lambda_1^hA_3+\lambda_2^hB_3)\right|^\alpha}\Gamma(ds)\\\sim\alpha D_4\,\lambda_1^h\ \ \text{\normalsize for}\ \ h \to +\infty,
						\end{multline}}	
					where 
					{\footnotesize{\begin{align} \label{D_4}
							D_4:=\sum_{j=0}^{+\infty}\int_{S_2} A_3C_3^{\langle\alpha-1\rangle}\Gamma(ds)<+\infty.
							\end{align}}}
					{Let us mention that to apply the dominated convergence theorem the following conditions have to hold}
					{\footnotesize  \begin{equation}\begin{aligned}\label{conditions5} 
						&\sum_{j=0}^{+\infty}\int_{S_2}\left|A_3\right|^\alpha\Gamma(ds)<+\infty, \quad \sum_{j=0}^{+\infty}\int_{S_2}\left|A_3\right||C_3|^{\alpha-1}\Gamma(ds)<+\infty,\\
						&\sum_{j=0}^{+\infty}\int_{S_2}\left|B_3\right|^{\alpha}\Gamma(ds)<+\infty, \quad
						\sum_{j=0}^{+\infty}\int_{S_2}\left|B_3\right||C_3|^{\alpha-1}\Gamma(ds)<+\infty,
						\end{aligned}\end{equation}}
					{which is always true (see Remark \ref{remark} in Appendix C).}
					\item Let us consider the case of $|\lambda_1|<|\lambda_2|$. Again, proceeding exactly as above we obtain
					{\footnotesize  \begin{multline}
						\sum_{j=0}^{+\infty}\int_{S_2}{\left|\lambda_1^hA_3+\lambda_2^hB_3\right|^{\alpha}+\left|C_3\right|^{\alpha}-\left|C_3-(\lambda_1^hA_3+\lambda_2^hB_3)\right|^\alpha}\Gamma(ds)\\\sim\alpha D_5\,\lambda_2^h\ \ \text{\normalsize for}\ \ h \to +\infty,
						\end{multline}}	
					where 
					{\footnotesize  \begin{align} \label{D_5} D_5:=\sum_{j=0}^{+\infty}\int_{S_2} B_3C_3^{\langle\alpha-1\rangle}\Gamma(ds)<+\infty. 
						\end{align}}
					\item Let us consider the case of $\lambda_1=-\lambda_2$ and $h$ is even. Proceeding exactly as above we have
					{\footnotesize  \begin{multline}
						\sum_{j=0}^{+\infty}\int_{S_2}{\left|\lambda_1^hA_3+\lambda_2^hB_3\right|^{\alpha}+\left|C_3\right|^{\alpha}-\left|C_3-(\lambda_1^hA_3+\lambda_2^hB_3)\right|^\alpha}\Gamma(ds)\\\sim\alpha (D_4+D_5)\,\lambda_1^h\ \text{\normalsize for}\ h \to +\infty,
						\end{multline}}	
					where 
					{\footnotesize  \begin{align*}
						D_4+D_5=\sum_{j=0}^{+\infty}\int_{S_2} (A_3+B_3)C_3^{\langle\alpha-1\rangle}\Gamma(ds)<+\infty. \end{align*}}
					\item Let us consider the case of $\lambda_1=-\lambda_2$ and $h$ is odd. Proceeding exactly as above we have
					{\footnotesize  \begin{multline}
						\sum_{j=0}^{+\infty}\int_{S_2}{\left|\lambda_1^hA_3+\lambda_2^hB_3\right|^{\alpha}+\left|C_3\right|^{\alpha}-\left|C_3-(\lambda_1^hA_3+\lambda_2^hB_3)\right|^\alpha}\Gamma(ds)\\\sim\alpha (D_4-D_5)\,\lambda_1^h\  \text{\normalsize for}\ h \to +\infty,
						\end{multline}}
					where 
					{\footnotesize  \begin{align*} D_4-D_5=\sum_{j=0}^{+\infty}\int_{S_2} (A_3-B_3)C_3^{\langle\alpha-1\rangle}\Gamma(ds)<+\infty.\end{align*}}
				\end{enumerate}
				\item Let us consider the case of real and equal eigenvalues of the coefficients' matrix $\Theta$, $\lambda_1=\lambda_2=\lambda$, and $|\lambda|<1$. For the cross-codifference function $\rm{CD}(X_1(t),X_2(t+h))$ given in Lemma \ref{th1} (see Eq. (\ref{eq9})) we take the notation
				{\footnotesize  \begin{equation}\begin{aligned}\label{A4B4C4}
					A_4(s_1,s_2,a_3,a_4,\lambda,j)=&A_4=j\lambda^{j-1}a_3s_1+j\lambda^{j-1}a_4s_2-j\lambda^js_2+\lambda^js_2,\\
					B_4(s_1,s_2,a_1,a_2,\lambda,j)=&B_4=\lambda^{j-1}a_3s_1+\lambda^{j-1}a_4s_2-\lambda^{j}s_2,\\
					C_4(s_1,s_2,a_3,a_4,\lambda,j)=&C_4=j\lambda^{j-1}a_1s_1-j\lambda^{j}s_1+\lambda^js_1+j\lambda^{j-1}a_2s_2.
					\end{aligned}\end{equation}}
				Referring to the above expressions as $A_4$, $B_4$ and $C_4$ to simplify the notation, we can write the cross-codifference function $\rm{CD}(X_1(t),X_2(t+h))$ given in Eq. (\ref{eq9}) in the following form
				{\footnotesize  \begin{multline}
					{\rm{CD}}(X_1(t),X_2(t+h))\\=\sum_{j=0}^{+\infty}\int_{S_2}\left(|	\lambda^h{A_4}+h\lambda^h{B_4}|^{\alpha}+|C_4|^{\alpha}-|C_4-(\lambda^h{A_4}+h\lambda^h{B_4})|^\alpha\right)\Gamma(ds).
					\end{multline}}
				Now, proceeding exactly as in a) we obtain
				{\footnotesize  \begin{multline}
					\sum_{j=0}^{+\infty}\int_{S_2} \left|\lambda^hA_4+h\lambda^hB_4\right|^{\alpha}+\left|C_4\right|^{\alpha}-\left|C_4-(\lambda^hA_4+h\lambda^hB_4)\right|^\alpha\Gamma(ds) \\ \sim \alpha D_6\,h\lambda^h \ \ \text{\normalsize for}\ \ h \to +\infty,
					\end{multline}}
				where 
				{\footnotesize  \begin{align} \label{D_6}
					D_6:=\sum_{j=0}^{+\infty}\int_{S_2} B_4C_4^{\langle\alpha-1\rangle}\Gamma(ds)<+\infty.
					\end{align}}
				{Similarly as above, to apply the dominated convergence theorem the following conditions have to hold}
				{\footnotesize  \begin{equation}
					\begin{aligned}\label{conditions7} 
					&\sum_{j=0}^{+\infty}\int_{S_2}\left|A_4\right|^\alpha\Gamma(ds)<+\infty, \quad \sum_{j=0}^{+\infty}\int_{S_2}\left|A_4\right||C_4|^{\alpha-1}\Gamma(ds)<+\infty,\\
					&\sum_{j=0}^{+\infty}\int_{S_2}\left|B_4\right|^{\alpha}\Gamma(ds)<+\infty, \quad
					\sum_{j=0}^{+\infty}\int_{S_2}\left|B_4\right||C_4|^{\alpha-1}\Gamma(ds)<+\infty,
					\end{aligned}
					\end{equation}}	
			{which is always satisfied (see Remark \ref{remark} in Appendix C).}
			\end{itemize}
	\end{enumerate}}

\section*{Appendix B: {Proof of Lemma \ref{lema2}}}

	{\begin{enumerate}[a)]
			\item First, we examine the asymptotic behavior of ${\rm{CV}}(X_1(t),X_2(t-h))$ for $h \to +\infty$. 
			\begin{itemize}
				\item Let us consider the case of two different real eigenvalues of the coefficients' matrix $\Theta$, $\lambda_1 \neq \lambda_2$, and $|\lambda_1|<1$, $|\lambda_2|<1$. By using the notation introduced in Eq. (\ref{A1B1C1}) the cross-covariation function ${\rm{CV}}(X_1(t),X_2(t-h))$ given in Lemma \ref{th2} (see Eq. (\ref{eq4})) can be written in the following form
				{\footnotesize  \begin{align} \label{CVminus}
					{\rm{CV}}(X_1(t),X_2(t-h))=\lambda_1^hD_1+\lambda_2^hD_2,
					\end{align}}
				where $D_1$ and $D_2$ are specified in Eqs. (\ref{D_1}) and (\ref{D_2}), respectively. The formula given in Eq. (\ref{CVminus}) directly leads to the expressions given in Lemma \ref{lema2}.
				\item Let us consider the case of the real and equal eigenvalues of the coefficients' matrix $\Theta$, $\lambda_1=\lambda_2=\lambda$, and $|\lambda|<1$. By using the notation introduced in Eq. (\ref{A2B2C2}) the cross-covariation function ${\rm{CV}}(X_1(t),X_2(t-h))$ given in Lemma \ref{th2} (see Eq. (\ref{eq6})) can be written in the following form
				{\footnotesize  \begin{align} \label{CVminus2}
					{\rm{CV}}(X_1(t),X_2(t-h))=\lambda^hE_3+h\lambda^hD_3,
					\end{align}}
				where $D_3$ is given in Eq. (\ref{D_3}) and
				{\footnotesize  \begin{align*}
					E_3={\sum_{j=0}^{+\infty}\int_{S_2}A_2C_2^{\langle \alpha -1 \rangle}\Gamma(ds)}.
					\end{align*}}
				Similarly to the previous case, the formula given in Eq. (\ref{CVminus2}) directly leads to the expression given in Lemma \ref{lema2}. 	
			\end{itemize}
			\item Now, we examine the asymptotic behavior of ${\rm{CV}}(X_1(t),X_2(t+h))$ for $h \to +\infty$.
			\begin{itemize}
				\item Let us consider the case of two different real eigenvalues of the coefficients' matrix $\Theta$, $\lambda_1 \neq \lambda_2$, and $|\lambda_1|<1$, $|\lambda_2|<1$. Taking the notation introduced in Eq. (\ref{A3B3C3}) the cross-covariation function ${\rm{CV}}(X_1(t),X_2(t+h))$ given in Lemma \ref{th2} (see Eq. (\ref{eq8})) can be written in the following form
				{\footnotesize  \begin{align}
					\rm{CV}(X_1(t),X_2(t+h))=\sum_{j=0}^{+\infty}\int_{S_2}C_3\left(A_3\lambda_1^h+B_3\lambda_2^h\right)^{\langle \alpha-1 \rangle}\Gamma(ds).
					\end{align}}
				\begin{enumerate}[I)]
					\item Let us consider the case of $|\lambda_1|>|\lambda_2|$. Now, similarly as in Lemma \ref{lema1}, we split the proof into two parts. Namely, we show that
					\begin{enumerate}[i)]
						\item 
						{\footnotesize  \begin{multline} \label{equality5}
							\lim_{h \to +\infty}\sum_{j=0}^{+\infty}\int_{S_2}\frac{C_3\left(A_3\lambda_1^h+B_3\lambda_2^h\right)^{\langle \alpha-1 \rangle}}{\left(\lambda_1^h\right)^{\langle\alpha-1\rangle}}\Gamma(ds)\\=\sum_{j=0}^{+\infty}\lim_{h \to +\infty}\int_{S_2}\frac{C_3\left(A_3\lambda_1^h+B_3\lambda_2^h\right)^{\langle \alpha-1 \rangle}}{\left(\lambda_1^h\right)^{\langle\alpha-1\rangle}}\Gamma(ds),
							\end{multline}}
						\item
						{\footnotesize  \begin{multline} \label{equality6}
							\lim_{h \to +\infty}\int_{S_2}\frac{C_3\left(A_3\lambda_1^h+B_3\lambda_2^h\right)^{\langle \alpha-1 \rangle}}{\left(\lambda_1^h\right)^{\langle\alpha-1\rangle}}\Gamma(ds)\\=\int_{S_2}\lim_{h \to +\infty}\frac{C_3\left(A_3\lambda_1^h+B_3\lambda_2^h\right)^{\langle \alpha-1 \rangle}}{\left(\lambda_1^h\right)^{\langle\alpha-1\rangle}}\Gamma(ds).
							\end{multline}}
					\end{enumerate}
					At first, let us focus on part i). According to the dominated convergence theorem, the equality given in Eq. (\ref{equality5}) is true if the infinite sum over $j$ converges uniformly. Let us notice that using the following inequality true for all $n,m \in \R, 1<\alpha<2$
					{\footnotesize  \begin{align} \label{formula2}
						\left|n+m\right|^{\alpha-1}\leq|n|^{\alpha-1}+|m|^{\alpha-1}
						\end{align}}
					one can show that for all $j \in \N$ we have
					{\footnotesize \begin{align*}
						&\left|\int_{S_2}\frac{C_3\left(A_3\lambda_1^h+B_3\lambda_2^h\right)^{\langle \alpha-1 \rangle}}{\left(\lambda_1^h\right)^{\langle\alpha-1\rangle}}\Gamma(ds)\right| \leq\int_{S_2}\left|\frac{C_3\left(A_3\lambda_1^h+B_3\lambda_2^h\right)^{\langle \alpha-1 \rangle}}{\left(\lambda_1^h\right)^{\langle\alpha-1\rangle}}\right|\Gamma(ds)\\\leq& \int_{S_2}\frac{|C_3|\left|A_3\lambda_1^h+B_3\lambda_2^h\right|^{ \alpha-1}}{|\lambda_1^h|^{\alpha-1}}\Gamma(ds) = \int_{S_2}\frac{|C_3||\lambda_1^h|^{\alpha-1}\left|A_3+B_3\left(\frac{\lambda_2}{\lambda_1}\right)^h\right|^{\alpha-1}}{|\lambda_1^h|^{\alpha-1}}\Gamma(ds)\\=& \int_{S_2}{|C_3|\left|A_3+B_3\left(\frac{\lambda_2}{\lambda_1}\right)^h\right|^{\alpha-1}}\Gamma(ds)\\\leq&\int_{S_2}{|C_3|\left(\left|A_3\right|^{\alpha-1}+\left|B_3\right|^{\alpha-1}\left|\left(\frac{\lambda_2}{\lambda_1}\right)^h\right|^{\alpha-1}\right)}\Gamma(ds)\\\leq&\int_{S_2}|C_3|\left|A_3\right|^{\alpha-1}\Gamma(ds)+\int_{S_2}|C_3|\left|B_3\right|^{\alpha-1}\Gamma(ds)=K_j
						\end{align*}}
					and $K_j$ is independent of $h$. Now, the infinite sum under the limit in Eq. (\ref{equality5}) converges uniformly if the infinite sum of $M_j$ over $j \in \N$ converges, which is equivalent to the following conditions
					{\footnotesize  
						\begin{align}\label{conditions9} 
						\sum_{j=0}^{+\infty}\int_{S_2}\left|C_3\right||A_3|^{\alpha-1}\Gamma(ds)<+\infty,
						\quad \sum_{j=0}^{+\infty}\int_{S_2}\left|C_3\right||B_3|^{\alpha-1}\Gamma(ds)<+\infty,
						\end{align}}
					{which are always satisfied (see Remark \ref{remark} in Appendix C)} and thus the equality given in Eq. (\ref{equality5}) is true. Let us focus on part ii). According to the dominated convergence theorem, the equality given in Eq. (\ref{equality6}) is satisfied if the integrand is dominated by an integrable function independent of $h$. Using the formula given in Eq. (\ref{formula2}) for a fixed $\mathbf{s}=(s_1,s_2)\in S_2$ we obtain 
					{\footnotesize \begin{align*}
						&\left|\frac{C_3\left(A_3\lambda_1^h+B_3\lambda_2^h\right)^{\langle \alpha-1 \rangle}}{\left(\lambda_1^h\right)^{\langle\alpha-1\rangle}}\right| \leq \frac{|C_3|\left|A_3\lambda_1^h+B_3\lambda_2^h\right|^{\alpha-1}}{\left|\lambda_1^h\right|^{\alpha-1}} \\ = & \frac{|C_3||\lambda_1^h|^{\alpha-1}\left|A_3+B_3\left(\frac{\lambda_2}{\lambda_1}\right)^h\right|^{\alpha-1}}{\left|\lambda_1^h\right|^{\alpha-1}} \leq |C_3|\left(|A_3|^{\alpha-1}+\left|B_3\right|^{\alpha-1}\left|\left(\frac{\lambda_2}{\lambda_1}\right)^h\right|^{\alpha-1}\right)\\ \leq & |C_3||A_3|^{\alpha-1}+|C_3||B_3|^{\alpha-1}
						\end{align*}}
					which is independent of $h$. Now, since the dominating function should be integrable we obtain the following conditions for all $j \in \N$
					{\footnotesize  
						\begin{equation*}
						\begin{aligned}
						\int_{S_2}\left|C_3||A_3\right|^{\alpha-1}\Gamma(ds)<+\infty, \quad \int_{S_2}\left|C_3||B_3\right|^{\alpha-1}\Gamma(ds)<+\infty.\\
						\end{aligned}
						\end{equation*}}
					{Let us notice that the above conditions are satisfied since the conditions given in Eq. (\ref{conditions9}) hold.} Now, we calculate the limit of the integrand given in Eq. (\ref{equality6}) for a fixed $\mathbf{s}=(s_1,s_2)\in S_2$, namely
					{\footnotesize  \begin{align*} 
						\lim_{h \to +\infty}\frac{C_3\left(A_3\lambda_1^h+B_3\lambda_2^h\right)^{\langle \alpha-1 \rangle}}{\left(\lambda_1^h\right)^{\langle\alpha-1\rangle}}=\lim_{h \to +\infty}{C_3\left(A_3+B_3\left(\frac{\lambda_2}{\lambda_1}\right)^h\right)^{\langle \alpha-1 \rangle}}=C_3A_3^{\langle \alpha-1 \rangle}.
						\end{align*}}
					Now, {since the equalities presented in i) and ii) are true}, and we obtain
					{\footnotesize  \begin{align}
						\lim_{h \to +\infty}\sum_{j=0}^{+\infty}\int_{S_2}\frac{C_3\left(A_3\lambda_1^h+B_3\lambda_2^h\right)^{\langle \alpha-1 \rangle}}{\left(\lambda_1^h\right)^{\langle\alpha-1\rangle}}\Gamma(ds)= D_7.
						\end{align}}
					which is equivalent to the fact that
					{\footnotesize  \begin{align}
						\sum_{j=0}^{+\infty}\int_{S_2}C_3\left(A_3\lambda_1^h+B_3\lambda_2^h\right)^{\langle \alpha-1 \rangle}\Gamma(ds)\sim D_7\,\left(\lambda_1^h\right)^{\langle \alpha-1 \rangle}\ \ \text{\normalsize for}\ \ h \to +\infty,
						\end{align}}
					where
					{\footnotesize  \begin{align} \label{D_7}
						D_7:=\sum_{j=0}^{+\infty}\int_{S_2} C_3A_3^{\langle\alpha-1\rangle}\Gamma(ds).
						\end{align}}
					Let us notice that the conditions given in Eq. (\ref{conditions9}) guarantee that $D_7<+\infty$.
					\item Let us consider the case of $|\lambda_1|<|\lambda_2|$. Proceeding in the same way as above, we have
					{\footnotesize  \begin{align}
						\sum_{j=0}^{+\infty}\int_{S_2}C_3\left(A_3\lambda_1^h+B_3\lambda_2^h\right)^{\langle \alpha-1 \rangle}\Gamma(ds)\sim D_8\,\left(\lambda_2^h\right)^{\langle\alpha-1\rangle}\ \ \text{\normalsize for}\ \ h \to +\infty,
						\end{align}}
					where
					{\footnotesize  \begin{align} \label{D_8}
						D_8:=\sum_{j=0}^{+\infty}\int_{S_2} C_3B_3^{\langle\alpha-1\rangle}\Gamma(ds) < +\infty.
						\end{align}} 
					\item Let us consider $\lambda_1=-\lambda_2$ and even $h$. In this case, we obtain that
					{\footnotesize  \begin{align}
						\sum_{j=0}^{+\infty}\int_{S_2}C_3\left(A_3\lambda_1^h+B_3\lambda_2^h\right)^{\langle \alpha-1 \rangle}\Gamma(ds) = D_{10}\,\left(\lambda_1^h\right)^{\langle \alpha-1\rangle},
						\end{align}}
					where
					{\footnotesize  \begin{align} \label{D_10}
						D_{10}=\sum_{j=0}^{+\infty}\int_{S_2}C_3(A_3+B_3)^{\langle\alpha-1\rangle}\Gamma(ds)<+\infty.\end{align}}
					\item Let us consider $\lambda_1=-\lambda_2$ and odd $h$. Similarly to the previous case, we have
					{\footnotesize  \begin{align}
						\sum_{j=0}^{+\infty}\int_{S_2}C_3\left(A_3\lambda_1^h+B_3\lambda_2^h\right)^{\langle \alpha-1 \rangle}\Gamma(ds)= D_{11}\,\left(\lambda_1^h\right)^{\langle \alpha-1\rangle},
						\end{align}}
					where
					{\footnotesize  \begin{align} \label{D_11}
						D_{11}=\sum_{j=0}^{+\infty}\int_{S_2}C_3(A_3-B_3)^{\langle\alpha-1\rangle}\Gamma(ds).
						\end{align}}
				\end{enumerate}
				\item Let us consider the case of the real and equal eigenvalues of the coefficients' matrix $\Theta$, $\lambda_1=\lambda_2=\lambda$, and $|\lambda|<1$. Using the notation introduced in Eq. (\ref{A4B4C4}) the cross-covariation function ${\rm{CV}}(X_1(t),X_2(t+h))$ given in Lemma \ref{th2} (see Eq. (\ref{eq10})) can be expressed as follows
				{\footnotesize  \begin{align}
					{\rm{CV}}(X_1(t),X_2(t+h))=\sum_{j=0}^{+\infty}\int_{S_2}C_4\left(A_4\lambda^h+B_4h\lambda^h\right)^{\langle \alpha-1 \rangle}\Gamma(ds).
					\end{align}}
				Now, similarly to the previous case we split the proof into two parts, namely we show that
				\begin{enumerate}[i)]
					\item 
					{\footnotesize  \begin{multline} \label{equality7}
						\lim_{h \to +\infty}\sum_{j=0}^{+\infty}\int_{S_2}\frac{C_4\left(A_4\lambda^h+B_4h\lambda^h\right)^{\langle \alpha-1 \rangle}}{\left(h\lambda^h\right)^{\langle\alpha-1\rangle}}\Gamma(ds)\\=\sum_{j=0}^{+\infty}\lim_{h \to +\infty}\int_{S_2}\frac{C_4\left(A_4\lambda^h+B_4h\lambda^h\right)^{\langle \alpha-1 \rangle}}{\left(h\lambda^h\right)^{\langle\alpha-1\rangle}}\Gamma(ds),
						\end{multline}}
					\item
					{\footnotesize  \begin{multline} \label{equality8}
						\lim_{h \to +\infty}\int_{S_2}\frac{C_4\left(A_4\lambda^h+B_4h\lambda^h\right)^{\langle \alpha-1 \rangle}}{\left(h\lambda^h\right)^{\langle\alpha-1\rangle}}\Gamma(ds)\\=\int_{S_2}\lim_{h \to +\infty}\frac{C_4\left(A_4\lambda^h+B_4h\lambda^h\right)^{\langle \alpha-1 \rangle}}{\left(h\lambda^h\right)^{\langle\alpha-1\rangle}}\Gamma(ds).
						\end{multline}}
				\end{enumerate}
				Now, according to the dominated convergence theorem, the equality corresponding to the part i), given in Eq. (\ref{equality7}), is satisfied if the infinite sum over $j \in \N$ converges uniformly. Using the inequality given in Eq. (\ref{formula2}) one can show that for all $j \in \N$ we have
				{\footnotesize \begin{align*}
					&\left|\int_{S_2}\frac{C_4\left(A_4\lambda^h+B_4h\lambda^h\right)^{\langle \alpha-1 \rangle}}{\left(h\lambda^h\right)^{\langle\alpha-1\rangle}}\Gamma(ds)\right| \leq\int_{S_2}\left|\frac{C_4\left(A_4\lambda^h+B_4h\lambda^h\right)^{\langle \alpha-1 \rangle}}{\left(h\lambda^h\right)^{\langle\alpha-1\rangle}}\right|\Gamma(ds)\\\leq& \int_{S_2}\frac{|C_4|\left|A_4\lambda^h+B_4h\lambda^h\right|^{ \alpha-1}}{|h\lambda^h|^{\alpha-1}}\Gamma(ds) = \int_{S_2}\frac{|C_4||h\lambda^h|^{\alpha-1}\left|\frac{A_4}{h}+B_4\right|^{\alpha-1}}{|h\lambda^h|^{\alpha-1}}\Gamma(ds)\\=& \int_{S_2}{|C_4|\left|\frac{A_4}{h}+B_4\right|^{\alpha-1}}\Gamma(ds)\leq\int_{S_2}{|C_4|\left(\left|\frac{A_4}{h}\right|^{\alpha-1}+\left|B_4\right|^{\alpha-1}\right)}\Gamma(ds)\\\leq&\int_{S_2}|C_4|\left|A_4\right|^{\alpha-1}\Gamma(ds)+\int_{S_2}|C_4|\left|B_4\right|^{\alpha-1}\Gamma(ds)=L_j
					\end{align*}}
				and $L_j$ is independent of $h$. Therefore, the infinite sum under the limit in Eq. (\ref{equality7}) converges uniformly if the infinite sum of $L_j$ over $j \in \N$ converges, which is equivalent to the following conditions
				{\footnotesize 
					\begin{align}\label{conditions11} 
					\sum_{j=0}^{+\infty}\int_{S_2}\left|C_4\right||A_4|^{\alpha-1}\Gamma(ds)<+\infty,
					\quad \sum_{j=0}^{+\infty}\int_{S_2}\left|C_4\right||B_4|^{\alpha-1}\Gamma(ds)<+\infty,
					\end{align}}
				{which are always satisfied (see Remark \ref{remark} in Appendix C).	Now, since the equality given in Eq. (\ref{equality7}) is true, we consider the part ii).} Similarly, according to the dominated convergence theorem, the equality given in Eq. (\ref{equality8}) is satisfied if the integrand is dominated by an integrable function independent of $h$. Again, by using the formula given in Eq. (\ref{formula2}) for a fixed $\mathbf{s}=(s_1,s_2)\in S_2$ we have 
				{\footnotesize \begin{align*}
					&\left|\frac{C_4\left(A_4\lambda^h+B_4h\lambda^h\right)^{\langle \alpha-1 \rangle}}{\left(h\lambda^h\right)^{\langle\alpha-1\rangle}}\right| \leq \frac{|C_4|\left|A_4\lambda^h+B_4h\lambda^h\right|^{\alpha-1}}{\left|h\lambda^h\right|^{\alpha-1}} \\ = & \frac{|C_4||h\lambda^h|^{\alpha-1}\left|\frac{A_4}{h}+B_4\right|^{\alpha-1}}{\left|h\lambda^h\right|^{\alpha-1}} \leq |C_4|\left(|\frac{A_4}{h}|^{\alpha-1}+\left|B_4\right|^{\alpha-1}\right)\\\leq& |C_4||A_4|^{\alpha-1}+|C_4||B_4|^{\alpha-1},
					\end{align*}}
				which is independent of $h$. Now, the fact that the dominating function should be integrable leads to the following conditions for all $j \in \N$
				{\footnotesize 
					\begin{equation*}
					\begin{aligned}
					\int_{S_2}\left|C_4||A_4\right|^{\alpha-1}\Gamma(ds)<+\infty, \quad \int_{S_2}\left|C_4||B_4\right|^{\alpha-1}\Gamma(ds)<+\infty,\\
					\end{aligned}
					\end{equation*}}
				{which are satisfied since the conditions given in (\ref{conditions11}) hold.} Now, we calculate the limit of the integrand given in Eq. (\ref{equality8}) for a fixed $\mathbf{s}=(s_1,s_2)\in S_2$, namely
				{\footnotesize \begin{align*} 
					\lim_{h \to +\infty}\frac{C_4\left(A_4\lambda^h+B_4h\lambda^h\right)^{\langle \alpha-1 \rangle}}{\left(h\lambda^h\right)^{\langle\alpha-1\rangle}}=\lim_{h \to +\infty}{C_4\left(\frac{A_4}{h}+B_4\right)^{\langle \alpha-1 \rangle}}=C_4B_4^{\langle \alpha-1 \rangle}.
					\end{align*}}
				Finally, {since the equalities presented in i) and ii) are true}, and we have
				{\footnotesize \begin{align}
					\lim_{h \to +\infty}\sum_{j=0}^{+\infty}\int_{S_2}\frac{C_4\left(A_4\lambda^h+B_4h\lambda^h\right)^{\langle \alpha-1 \rangle}}{\left(h\lambda^h\right)^{\langle\alpha-1\rangle}}\Gamma(ds)= D_9,
					\end{align}}
				which is equivalent to the fact that
				{\footnotesize \begin{align}
					\sum_{j=0}^{+\infty}\int_{S_2}C_4\left(A_4\lambda^h+B_4h\lambda^h\right)^{\langle \alpha-1 \rangle}\Gamma(ds)\sim D_9\,\left(h\lambda^h\right)^{\langle \alpha-1 \rangle}\ \ \text{\normalsize for}\ \ h \to +\infty,
					\end{align}}
				where
				{\footnotesize \begin{align} \label{D_9}
					D_9:=\sum_{j=0}^{+\infty}\int_{S_2} C_4B_4^{\langle\alpha-1\rangle}\Gamma(ds).
					\end{align}}
				Let us notice that the conditions given in Eq. (\ref{conditions11}) guarantee that $D_9<+\infty$.
			\end{itemize}
	\end{enumerate}}
{\section*{Appendix C}
	\begin{remark} \label{remark}
		Let us notice that since the moduli of $A_i, B_i, C_i$ for $i=1,2,3,4$ can be upper-bounded either by $M\max\left(|\lambda_1|,|\lambda_2|\right)^j$ or by $M\max j\left(|\lambda_1|,|\lambda_2|\right)^j$, where the constant $M$ is independent of $j$ and $\max \left(|\lambda_1|,|\lambda_2|\right)<1$, and since the spectral measure is finite, i.e. $\Gamma(S_2)<+\infty$, all conditions given in the proofs of Lemmas \ref{lema1} and \ref{lema2}, see Eqs. (\ref{conditions1}), (\ref{conditions5}), (\ref{conditions9}), (\ref{conditions3}), (\ref{conditions7})  and (\ref{conditions11}), are satisified for any spectral measure $\Gamma(\cdot)$.
\end{remark}}

\section*{Acknowledgments}
We would like to acknowledge the support of the National Center of Science Opus Grant No. 2016/21/B/ST1/00929 "Anomalous diffusion processes and their applications in real data modelling".

\bibliographystyle{plainnat}
\bibliography{bib}

\begin{thebibliography}{39}
\providecommand{\natexlab}[1]{#1}
\providecommand{\url}[1]{\texttt{#1}}
\expandafter\ifx\csname urlstyle\endcsname\relax
  \providecommand{\doi}[1]{doi: #1}\else
  \providecommand{\doi}{doi: \begingroup \urlstyle{rm}\Url}\fi

\bibitem[Anderson and Meerschaert(1998)]{infinite11}
Paul~L. Anderson and Mark~M. Meerschaert.
\newblock Modeling river flows with heavy tails.
\newblock \emph{Water Resources Research}, 34\penalty0 (9):\penalty0
  2271--2280, 1998.

\bibitem[Burnecki et~al.(2015)Burnecki, Wy{\l}oma{\'n}ska, and Chechkin]{phy2}
K.~Burnecki, A.~Wy{\l}oma{\'n}ska, and A.~Chechkin.
\newblock Discriminating between light- and heavy-tailed distributions with
  limit theorem.
\newblock \emph{PLoS ONE}, 10(12):\penalty0 1--23, 2015.

\bibitem[Butler and Okada(2009)]{aw2}
Kirt~C. Butler and Katsushi Okada.
\newblock The relative contribution of conditional mean and volatility in
  bivariate returns to international stock market indices.
\newblock \emph{Applied Financial Economics}, 19\penalty0 (1):\penalty0 1--15,
  2009.

\bibitem[Chen et~al.(2016)Chen, Geng, and Yin]{floc1}
Z.~Chen, X.~Geng, and F.~Yin.
\newblock A harmonic suppression method based on fractional lower order
  statistics for power system.
\newblock \emph{Industrial Electronics IEEE Transactions on}, 63(6):\penalty0
  3745--3755, 2016.

\bibitem[Damarackas and Paulauskas(2016)]{spectral_cov}
Julius Damarackas and Vygantas Paulauskas.
\newblock Spectral covariance and limit theorems for random fields with
  infinite variance.
\newblock \emph{Journal of Multivariate Analysis}, 153, 2016.

\bibitem[Gallagher(2001)]{est2}
C.~M. Gallagher.
\newblock A method for fitting stable autoregressive models using the
  autocovariation function.
\newblock \emph{Statistics \& Probability Letter}, 53:\penalty0 381--390, 2001.

\bibitem[Gardes and Girard(2010)]{Gardes2010}
Laurent Gardes and St{\'e}phane Girard.
\newblock Conditional extremes from heavy-tailed distributions: an application
  to the estimation of extreme rainfall return levels.
\newblock \emph{Extremes}, 13\penalty0 (2):\penalty0 177--204, 2010.

\bibitem[Glasserman et~al.(2002)Glasserman, Heidelberger, and
  Shahabuddin]{valueatrisk}
Paul Glasserman, Philip Heidelberger, and Perwez Shahabuddin.
\newblock Portfolio {V}alue-at-{R}isk with {H}eavy-{T}ailed {R}isk {F}actors.
\newblock \emph{Mathematical Finance}, 12\penalty0 (3):\penalty0 239--269,
  2002.

\bibitem[Grzesiek et~al.(2019{\natexlab{a}})Grzesiek, Teuerle, Sikora, and
  Wy{\l}oma{\'n}ska]{nasza2}
A.~Grzesiek, M.~Teuerle, G.~Sikora, and A.~Wy{\l}oma{\'n}ska.
\newblock Spatial-temporal dependence measures for $\alpha-$stable bivariate
  {AR}(1).
\newblock \emph{submitted in 2019}, 2019{\natexlab{a}}.

\bibitem[Grzesiek et~al.(2019{\natexlab{b}})Grzesiek, Teuerle, and
  Wy{\l}oma{\'n}ska]{nasza}
A.~Grzesiek, M.~Teuerle, and A.~Wy{\l}oma{\'n}ska.
\newblock Cross-codifference for bidimensional {VAR}(1) models with infinite
  variance.
\newblock 2019{\natexlab{b}}.
\newblock \ 27. Located at: arXiv:1902.02142.

\bibitem[Jab{\l}o{\'n}ska-Sabuka et~al.(2017)Jab{\l}o{\'n}ska-Sabuka, Teuerle,
  and Wy{\l}oma{\'n}ska]{stab_jablonska}
Matylda Jab{\l}o{\'n}ska-Sabuka, Marek Teuerle, and Agnieszka
  Wy{\l}oma{\'n}ska.
\newblock Bivariate sub-{G}aussian model for stock index returns.
\newblock \emph{Physica A: Statistical Mechanics and its Applications},
  486:\penalty0 628--637, 2017.

\bibitem[Kozubowski and Panorska(1999)]{aw3}
T.~J. Kozubowski and A.~K. Panorska.
\newblock Multivariate geometric stable distributions in financial
  applications.
\newblock \emph{Mathematical and Computer Modelling}, 29\penalty0
  (10-12):\penalty0 83--92, 1999.

\bibitem[Kozubowski et~al.(2003)Kozubowski, Panorska, and Rachev]{panor1}
T.~J. Kozubowski, A.~K. Panorska, and S.~T. Rachev.
\newblock Statistical issues in modeling multivariate stable portfolios.
\newblock In Svetlozar~T. Rachev, editor, \emph{Handbook of Heavy Tailed
  Distributions in Finance}, volume~1 of \emph{Handbooks in Finance}, pages 131
  -- 167. Amsterdam: North-Holland, 2003.

\bibitem[Kysel{\'y}(2010)]{Kysely2010}
Jan Kysel{\'y}.
\newblock Coverage probability of bootstrap confidence intervals in
  heavy-tailed frequency models, with application to precipitation data.
\newblock \emph{Theoretical and Applied Climatology}, 101\penalty0
  (3):\penalty0 345--361, 2010.

\bibitem[Lii and Rosenblatt(1992)]{multiarmastable}
K.-S. Lii and M.~Rosenblatt.
\newblock An approximate maximum likelihood estimation for non-gaussian
  non-minimum phase moving average processes.
\newblock \emph{Journal of Multivariate Analysis}, 43\penalty0 (2):\penalty0
  272--299, 1992.

\bibitem[Ma and L.~Nikias(1996)]{floc2}
X.~Ma and C.~L.~Nikias.
\newblock Joint estimation of time delay and frequency delay in impulsive noise
  using fractional lower order statistics.
\newblock \emph{Signal Processing, IEEE Transactions on}, 44(11):\penalty0
  2669--2687, 1996.

\bibitem[Maejima and Yamamoto(2003)]{maejima2003}
Makoto Maejima and Kenji Yamamoto.
\newblock Long-{M}emory {S}table {O}rnstein-{U}hlenbeck {P}rocesses.
\newblock \emph{Electronic Journal of Probability}, 8\penalty0 (19):\penalty0
  1--18, 2003.

\bibitem[Mittnik and Rachev(1991)]{stab4}
S.~Mittnik and S.~T. Rachev.
\newblock Alternative multivariate stable distributions and their applications
  to financial modeling.
\newblock In Cambanis S., Samorodnitsky G., and Taqqu M.S., editors,
  \emph{Stable Processes and Related Topics}, volume~25 of \emph{Progress in
  Probabilty}. Birkh{\"a}user Boston, 1991.

\bibitem[Mittnik and Rachev(2000)]{mit}
S.~Mittnik and S.~T. Rachev.
\newblock \emph{Stable Paretian Models in Finance}.
\newblock New York: Wiley, 2000.

\bibitem[Nolan and Panorska(1997)]{panor2}
J.~P. Nolan and A.~K. Panorska.
\newblock Data analysis for heavy tailed multivariate samples.
\newblock \emph{Communications in Statistics. Stochastic Models}, 13\penalty0
  (4):\penalty0 687--702, 1997.

\bibitem[Nowicka(1997)]{arma}
Joanna Nowicka.
\newblock Asymptotic behavior of the covariation and the codifference for
  {ARMA} models with stable innovations.
\newblock \emph{Communications in Statistics. Stochastic Models}, 13\penalty0
  (4):\penalty0 673--685, 1997.

\bibitem[Nowicka and Weron(1997)]{ARMA_measures}
Joanna Nowicka and Aleksander Weron.
\newblock Measures of dependence for {ARMA} models with stable innovations.
\newblock \emph{Annales Universitatis Mariae Curie-Skłodowska. Sectio A.
  Mathematica}, LI 1,14:\penalty0 133--144, 1997.

\bibitem[Nowicka and Wy{\l}oma{\'n}ska(2006)]{parma}
Joanna Nowicka and Agnieszka Wy{\l}oma{\'n}ska.
\newblock The dependence structure for {PARMA} models with a-stable
  innovations.
\newblock \emph{Acta Physica Polonica B}, 37\penalty0 (11):\penalty0
  3071--3081, 2006.

\bibitem[Nowicka-Zagrajek and Wy{\l}oma{\'n}ska(2008)]{ar}
J.~Nowicka-Zagrajek and A.~Wy{\l}oma{\'n}ska.
\newblock Measures of dependence for stable {AR}(1) models with time-varying
  coefficients.
\newblock \emph{Stochastic Models}, 24\penalty0 (1):\penalty0 58--70, 2008.

\bibitem[Obuchowski and A.(2013)]{fin2}
J.~Obuchowski and Wy{\l}oma{\'n}ska A.
\newblock The {O}rnstein-{U}hlenbeck process with non-{G}aussian structure.
\newblock \emph{Acta Physica Polonica B}, 44(5):\penalty0 1123--1136, 2013.

\bibitem[Press(1972)]{stab3}
S.J. Press.
\newblock Multivariate stable distributions.
\newblock \emph{Journal of Multivariate Analysis}, 2\penalty0 (4):\penalty0
  444--462, 1972.

\bibitem[Rosadi(2009)]{cod2}
D.~Rosadi.
\newblock Testing for independence in heavy-tailed time series using the
  codifference function.
\newblock \emph{Computational Statistics \& Data Analysis}, 53(12):\penalty0
  4516--4529, 2009.

\bibitem[Rosadi and M.(2011)]{cod1}
D.~Rosadi and Deistler M.
\newblock Estimating the codifference function of linear time series models
  with infinite variance.
\newblock \emph{Metrika}, 73(3):\penalty0 395--429, 2011.

\bibitem[Rosadi(2016)]{dep7}
Dedi Rosadi.
\newblock Measuring dependence of random variables with finite and infinite
  variance using the codifference and the generalized codifference function.
\newblock \emph{AIP Conference Proceedings}, 1755\penalty0 (1):\penalty0
  120004, 2016.

\bibitem[Saichev and Sornette(2004)]{earthquake}
A.~Saichev and D.~Sornette.
\newblock Anomalous power law distribution of total lifetimes of branching
  processes: Application to earthquake aftershock sequences.
\newblock \emph{Physical Review E}, 70:\penalty0 046123, 2004.

\bibitem[Samorodnitsky and Taqqu(1994)]{Taqqu}
G.~Samorodnitsky and M.~S. Taqqu.
\newblock \emph{Stable Non-Gaussian Random Processes: Stochastic Models with
  Infinite Variance}.
\newblock New York: Chapman \& Hall, 1994.

\bibitem[Shao and Nikias(1993)]{floc3}
M.~Shao and C.~L. Nikias.
\newblock Signal processing with fractional lower order moments: stable
  processes and their applications.
\newblock \emph{Proceedings of the IEEE}, 81\penalty0 (7):\penalty0 986--1010,
  1993.

\bibitem[Stoyanov et~al.(2006)Stoyanov, Samorodnitsky, Rachev, and
  Ortobelli]{stab_stoyanov}
Stoyan~V. Stoyanov, Gennady Samorodnitsky, Svetlozar Rachev, and Sergio
  Ortobelli.
\newblock Computing the {p}ortfolio {c}onditional {V}alue-at-{R}isk in the
  alpha-stable case.
\newblock \emph{Probability and Mathematical Statistics}, 26:\penalty0 1--22,
  2006.

\bibitem[Weir(1973)]{dom_con_theorem}
A.~J. Weir.
\newblock \emph{Lebesgue Integration and Measure}.
\newblock Cambridge: Cambridge University Press, 1973.

\bibitem[Williams(1992)]{matrix}
K.~S. Williams.
\newblock The nth power of a $2 \times 2$ matrix.
\newblock \emph{Mathematics Magazine}, 65(5):\penalty0 336--336, 1992.

\bibitem[Wy{\l}oma{\'n}ska et~al.(2015)Wy{\l}oma{\'n}ska, Chechkin, Sokolov,
  and Gajda]{Chechkin2015}
A.~Wy{\l}oma{\'n}ska, A.~Chechkin, I.M. Sokolov, and J.~Gajda.
\newblock Codifference as a practical tool to measure interdependence.
\newblock \emph{Physica A}, 421:\penalty0 412--429, 2015.

\bibitem[{\.Z}ak et~al.(2017){\.Z}ak, Wy{\l}oma{\'n}ska, and R.]{dia2}
G.~{\.Z}ak, A.~Wy{\l}oma{\'n}ska, and Zimroz R.
\newblock Data driven iterative vibration signal enhancement strategy using
  alpha stable distribution.
\newblock \emph{Shock and Vibration}, 2017:\penalty0 11, 2017.

\bibitem[{\.Z}ak et~al.(2014){\.Z}ak, Obuchowski, Wy{\l}oma{\'n}ska, and
  Zimroz]{stab_obuchowski}
Grzegorz {\.Z}ak, Jakub Obuchowski, Agnieszka Wy{\l}oma{\'n}ska, and
  Rados{\l}aw Zimroz.
\newblock Application of {ARMA} modelling and alpha-stable distribution for
  local damage detection in bearings.
\newblock \emph{Diagnostyka}, 15, 01 2014.

\bibitem[Zolotarev(1986)]{zolotarev}
V.~M. Zolotarev.
\newblock \emph{One-dimensional stable distributions}.
\newblock Translations of Mathematical Monographs. Providence: American
  Mathematical Society, 1986.

\end{thebibliography}

\end{document}